\documentclass[preprint,12pt]{elsarticle}

\usepackage[top=2.5cm,bottom=2.5cm,left=2.0cm,right=2.0cm]{geometry}

\usepackage{amssymb}
\usepackage{amsmath}
\usepackage{float}
\usepackage{tikz}
\usetikzlibrary{positioning}
\usepackage{verbatim}
\usetikzlibrary{arrows.meta,
                calc, chains,
                decorations.pathreplacing,
                calligraphy,
                positioning}
\usepackage[utf8]{inputenc}
\usepackage{bm}
\newtheorem{remark}{Remark}[section]

\def\ljump{{[\![}}
\def\rjump{{]\!]}}

\journal{}

\begin{document}

\begin{frontmatter}



\title{A Mesh-free Method Using Piecewise Deep Neural Network for Elliptic Interface Problems}
\author{Cuiyu He\fnref{label2}}
\ead{cuiyu.he@uga.edu}
\author{Xiaozhe Hu\fnref{label1}}
\ead{Xiaozhe.Hu@tufts.edu}
\author{Lin Mu\fnref{label2}}
\ead{linmu@uga.edu}
\fntext[label2]{Department of Mathematics, University of Georgia, Athens, GA 30602}
\fntext[label1]{Department of Mathematics, Tufts University, Medford, MA 02155}


%

\begin{abstract}
In this paper, we propose a novel mesh-free numerical method for solving the elliptic interface problems based on deep learning.  We approximate the solution by the neural networks and, since the solution may change dramatically across the interface, we employ different neural networks in different sub-domains. By reformulating the interface problem as a least-squares problem, we discretize the objective function using mean squared error via sampling and solve the proposed deep least-squares method by standard training algorithms such as stochastic gradient descent. The discretized objective function utilizes only the point-wise information on the sampling points and thus no underlying mesh is required. Doing this circumvents the challenging meshing procedure as well as the numerical integration on the complex interface.
To improve the computational efficiency for more challenging problems,
we further design an adaptive sampling strategy based on the residual of the least-squares function and propose an adaptive algorithm.  Finally, we present several numerical experiments in both 2D and 3D to show the flexibility, effectiveness, and accuracy of the proposed deep least-square method for solving interface problems.
\end{abstract}




\end{frontmatter}


\section{Introduction}
\label{Sect:Introduction}
Partial differential equations (PDEs) model a variety of physical and chemical phenomena such as diffusion, electrostatics, heat transfer, fluid dynamics, elasticity, and multi-phase flow in porous media.  Due to the complex nature of the PDEs, numerical simulations are oftentimes the only possible way for scientific discovery.  During the past century, many numerical methods have been developed to numerically solve the PDEs, for example, finite difference method \cite{leveque2007finite}, finite element method \cite{brenner2007mathematical}, finite volume method \cite{versteeg2007introduction}, spectral method \cite{xiu2010numerical}, and mesh-free method \cite{liu2009meshfree}. Moreover, the numerical solutions for PDEs have been widely used in many application fields, e.g., biology, petroleum engineering, meteorology, etc., and have achieved great success in the past several decades. 

However, there are still many challenging PDEs require the development of advanced numerical methods, e.g., turbulence flow, high-dimensional PDEs, interface problems, etc..  In this paper, we focus on the second-order elliptic interface problem, which captures many fundamental physical phenomena \cite{HANSBO20025537,doi:10.1137/0731054,doi:10.1137/130912700,MU2013106,peskin_2002,ZHOU20061}. There are mainly two categories of numerical methods, i.e., interface-fitted and -unfitted approaches. The first type uses interface-fitted meshes and such approach includes classical finite element and finite volume methods \cite{babuvska1970finite,bramble1996finite}, discontinuous Galerkin method \cite{massjung2012unfitted,cai2011discontinuous}, virtual element methods \cite{chen2017interface}.  The second type uses interface-unfitted meshes, for example, structured uniform meshes. Such an approach draws increasing attention during the past decade because it is difficult, if not impossible, to generate interface-fitted meshes, especially when the interface is geometrically complicated and/or time-dependent.  Typical methods of this type include immersed boundary methods \cite{peskin_2002}, immersed interface methods \cite{doi:10.1137/0731054}, matched interface and boundary method \cite{ZHOU20061}, ghost fluid method \cite{liu2000boundary}, extended finite element method \cite{fries2010extended}, cut finite element methods \cite{burman2015cutfem}, multi-scale FEM \cite{efendiev2009multiscale}, and immersed FEM \cite{doi:10.1137/130912700}.  Although both types of methods have been successful for solving interface problems to a certain extent, the implementation of those numerical schemes is not a straightforward task due to the jump conditions on the interface. In practice, interface problems remain quite hard due to the complicated geometry of the interfaces, which oftentimes are dynamically changing, and the singularities introduced by the interface conditions. Furthermore, neither type of method yields satisfactory results for the high dimensional interface problems, the non-linear interface problems, and other general interface problems. 

On the other hand, the neural network models have shown remarkable success in computer vision~\cite{handa2016gvnn}, pattern recognition~\cite{pao1989adaptive}, natural language processing \cite{collobert2008unified}, and many other artificial intelligence tasks. Despite being an old idea, the deep neural network (DNN) model also has great potential in nonlinear approximation, especially in modeling complicated data sets.  
The astonishing success of the DNN models in machine learning encourages wide applications to other fields, including recent studies of using the DNN models to numerically solve PDEs, especially those challenging ones which cannot be handled by existing numerical methods robustly and efficiently, e,g., \cite{CaiChenLiuLiu2019,DissanayakePhan-Thien1994,EYu,HanJentzenE,HeLiXuZheng,SamaniegoGoswami,SirignanoSpilopoulos,TranHamiltonMckayQuiringVassilevski}. 
The approximation properties of the DNN models, however, remain an active and open question.
Mathematically, there is a universal approximation theory about the single-layer neural network (see \cite{Pinkus} and references therein), which leads to more recent work~\cite{Zhou,SigelXu,DaubechiesDeVoreFoucartHanninPetrova}. In this paper, we focus on the numerical algorithm development and consider the theoretical approximation analysis as future work.


In this work, we use deep learning methods to solve interface problems. Our work is based on the recent work~\cite{EYu,CaiChenLiuLiu2019,WangZhang} by rewriting the second-order elliptic interface problem as a minimization problem.  More precisely, we use the least-squares (LS) approach to reformulate the interface problem and use the DNN model to approximate the solutions.  However, instead of using only one DNN structure to represent the numerical solution on the whole domain, we use two DNN structures to approximate the solution when the interface divides the domain into two sub-domains. 
This idea is based on the observation that the solution could undergo large jumps in derivative(s) across the interface and, therefore,
using one DNN structure could be inefficient to capture the difference.
 The numerical test results show that the proposed method is able to provide satisfactory 
 approximations of the solutions that  even have singularities on the interface. 
 Our approach can be considered as a piece-wise approximation and can be easily extended to more complicated interface problems with multiple sub-domains in which case we can use a piece-wise DNN structure in each sub-domain.  To solve the LS problem using the piece-wise DNN, we firstly sample some points and then define a discrete LS problem which can be solved by the stochastic gradient descent (SGD) method \cite{RobbinsMonro1951a}. 
 The advantage of using our discrete LS problem circumvents the meshing procedure that remains as a challenging task for problems with complex interfaces. Furthermore, on the interface, only the location information of the sampling points are required in the discretized formulation, the numerical integration on the entire interface is therefore alleviated. For more challenging problems with singularities, we further design an adaptive procedure that selects sampling points based on the point-wise value of the LS residual function. Numerical results have shown great improvement comparing to the uniform sampling strategy. 

The rest of the paper is organized as follows.  The second-order elliptic interface problem and its least-squares formulation are discussed in Section~\ref{Sect:Formulation}.  In Section~\ref{sec:deep-LS}, we introduce our deep least-squares method for solving the interface problem in detail.  Numerical results are shown in Section~\ref{sec:numerics} to demonstrate the efficiency of the proposed method.  Finally, we give some conclusions in Section~\ref{sec:conclusions}.

\section{Problem Formulation}
\label{Sect:Formulation}
In this section, we introduce the model interface problem as well as its classical LS formulation. 
For the sake of simplicity, we focus on the case that there is only one closed interface and the domain is divide into two sub-domains.  However, as noted our approach can be easily extended to more general cases.

\subsection{Interface Problem} \label{sec:interface-problem}
Let $\Omega$ be a bounded domain in $\mathbb{R}^d$, $d=2,3$, with Lipschitz boundary $\partial\Omega$, and the interface $\Gamma$ 
is closed and divides $\Omega$ into two disjoint sub-domains $\Omega_1$ and $\Omega_2$. We assume the interface is Lipschitz, however, our approach can handle more general interfaces in the same fashion. We consider the following second-order scalar elliptic interface problem,
\begin{eqnarray}
-\nabla\cdot(\beta(\bm{x})\nabla u) &=& f,\quad \mbox{ in }\Omega_1\cup\Omega_2,\label{eq:pde}\\
\ljump u \rjump &=& g_j,\quad \mbox{ on }\Gamma, \label{eqn:ic-u}\\
\ljump \beta(\bm{x})\nabla u\cdot\bm{n}\rjump &=& g_f,\quad \mbox{ on }\Gamma, \label{eqn:ic-flux}\\
u&=& g_D, \quad \mbox{ on }\partial\Omega,\label{eq:pde-bc}
\end{eqnarray}
where $f\in L^2(\Omega)$, $g_D \in H^{1/2}(\partial\Omega)$, and $\bm{n}$ is the unit outer normal vector to the interface $\Gamma$. The diffusion coefficient $\beta(\bm{x})  \ge \beta_0 > 0$ is a piece-wise constant function, i.e.,
\begin{eqnarray*}
\beta(\bm{x}) = \begin{cases}
\beta_1,\quad \mbox{ if } \bm{x} \in\Omega_1,\\
\beta_2,\quad \mbox{ if } \bm{x} \in\Omega_2,
\end{cases}
\end{eqnarray*} 
which has a finite jump of function value across the interface $\Gamma$. The notation $\ljump \cdot \rjump$ denotes the jump across the interface $\Gamma$ and is defined as follows,
\begin{equation*}
\ljump  u(\bm{x}) \rjump|_{\Gamma} = u_1(\bm{x})|_{\Gamma} - u_2(\bm{x})|_{\Gamma}, 
\end{equation*}
where 
\begin{equation*}
u(\bm{x}) = 
\begin{cases}
u_1(\bm{x}), \quad \mbox{ if } \bm{x} \in \Omega_1, \\
u_2(\bm{x}), \quad \mbox{ if } \bm{x} \in \Omega_2.
\end{cases}
\end{equation*}
For convenience, we refer $\Omega_1$ and $\Omega_2$ as the interior and exterior sub-domains, respectively.

\subsection{The Classical Least-squares Formulation} \label{sec:LS-formulation}
Throughout this paper, we shall use the standard notation and definitions for the Sobolev spaces $H^s(\Omega)$ and $H^s(\partial\Omega)$. The standard associated inner products are denoted by $(\cdot,\cdot)_{s,\Omega}$ and $(\cdot,\cdot)_{s,\Gamma}$ in $\Omega \in R^d$ and on $\Gamma \in R^{d-1}$, respectively.  And the standard induced norms are denoted by $\|\cdot\|_{s,\Omega}$ and $\|\cdot\|_{s,\Gamma}$. When $s = 0$, $H^0(\Omega)$ coincides with $L^2(\Omega)$.  When there is no ambiguity, the subscript $\Omega$ in the designation of norms will be suppressed.

To take advantage of the deep neural network, it is natural to consider LS formulation.  There are many different LS formulations for the interface problem~\eqref{eq:pde}-\eqref{eq:pde-bc}. One approach is to use the underlying minimization principle for elliptic interface problems. However, such an approach is limited to problems that have underlying minimization principle.  In this work, we adopt the simple LS principle proposed in~\cite{DissanayakePhan-Thien1994} and propose a LS functional that incorporates the interface conditions~\eqref{eqn:ic-u} and~\eqref{eqn:ic-flux} and the boundary condition~\eqref{eq:pde-bc} naturally. The LS functional is defined as follows,
\begin{eqnarray}
\mathcal{J}(v;g_j,g_f,g_D,f) &=& \|-\nabla\cdot\beta(\bm{x})\nabla v-f\|_{0,\Omega}^2+ {\beta}_j\|\ljump v\rjump-g_j\|_{0,\Gamma}^2\notag\\
&&+{\beta_f}\|\ljump\beta(\bm{x})\nabla v\cdot\bm{n} \rjump-g_f\|_{0,\Gamma}^2+\alpha\|v-g_D\|_{0,\partial\Omega}^2, \label{eqn:LS-obj}
\end{eqnarray}
for all $v\in H^1(\Omega)$, where {$\beta_j, \beta_f$ and $\alpha$ are constants} to be determined and, for the sake of simplicity, it may be chosen to be one.  The corresponding LS solution is then to find $u\in H^1(\Omega)$ such that
\begin{eqnarray}\label{eqn:LS-problem}
\mathcal{J}(u;g_j,g_f,g_D,f) = \min_{v\in H^1(\Omega)}\mathcal{J}(v;g_j,g_f,g_D,f).
\end{eqnarray}

\begin{remark}
By the Sobolev trace theorem, the interior and boundary terms in the functional $\mathcal{J}$ are not on the same scale.  As suggested in~\cite{CaiChenLiuLiu2019}, $\| \cdot \|_{3/2, \partial \Omega}$ can be used on the boundary to obtain a balanced LS functional.  Similarly, $\| \cdot \|_{3/2, \Gamma}$ can be used for the interface term.  In this work, we simply use the $L^2$ norm for both the boundary and interface terms and remark that our numeric experiments produce similar results.
\end{remark}

\section{Mesh-free Method Using DNN for interface problems}
\label{sec:deep-LS}
In this section, we discuss how to use DNN to numerically solve the interface problem based on the LS formulation~\eqref{eqn:LS-problem}.  The main idea of our new method is to use two neural networks to approximate the solution on two sub-domains, i.e., $\Omega_1$ and $\Omega_2$. This allows us to handle complicated interface problems by only looking at the location of the sampling points without the need for an underlying mesh.

\subsection{Deep Neural Network Structure}\label{Sect:DNN}
We first discuss the deep neural network structure used to approximate the solution $u(\bm{x})$.  A DNN structure is the composition of multiple linear functions and nonlinear activation functions. Specifically, the first component of DNN is a linear transformation $\bm{T}^\ell:\mathbb{R}^{n_\ell}\to\mathbb{R}^{n_{\ell+1}}$, {$\ell=1, \cdots, L$}, defined as follows,
\begin{eqnarray*}
	\bm{T}^\ell(\bm{x}^\ell) = \bm{W}^\ell \bm{x}^{\ell} +  \bm{b}^\ell,\mbox{ for }\bm{x}^\ell\in \mathbb{R}^{n_\ell},
\end{eqnarray*}
where $\bm{W}^\ell=(w_{i,j}^\ell)\in\mathbb{R}^{n_{\ell+1} \times n_\ell}$ and $\bm{b}^{\ell}\in\mathbb{R}^{n_{\ell+1}}$ are parameters in the DNN. The second component is an activation function $\psi:\mathbb{R}\to\mathbb{R}$ to be chosen.  Typical examples of the activation functions are $\tanh$, Sigmoid, and ReLU. Application of $\psi$ to a vector $\bm{x}\in\mathbb{R}^n$ is defined component-wisely, i.e., $\psi(\bm{x})=(\psi(x_i))$, $i=1,2,\cdots,n$. Then, the $\ell$-th layer of the DNN can be represented as the composition of the linear transform $\bm{T}^\ell$ and the nonlinear activation function $\psi$, i.e., 
\begin{equation*}
\mathcal{N}^\ell(\bm{x}^\ell) = \psi(\bm{T}^\ell(\bm{x}^\ell)),
\quad l=1, \cdots, L.
\end{equation*}
Note $\mathcal{N}^{\ell}: \mathbb{R}^{n_\ell} \mapsto \mathbb{R}^{n_{\ell+1}}$. A $L$-layer DNN is then defined as the composition of all $\mathcal{N}^{\ell}$, $\ell=1,2,\cdots,L$.  In particular, for an input $\bm{x}\in\mathbb{R}^{n_1}$, a general $L$-layer DNN can be represented as follows,
\begin{equation}\label{eqn:DNN}
\mathcal{NN}(\bm{x}; \Theta) = \bm{T}^L \circ \mathcal{N}^{L-1} \circ \cdots \circ \mathcal{N}^{2} \circ \mathcal{N}^1(\bm{x}),
\end{equation}
where $\Theta \in \mathbb{R}^N$ stands for all the parameters in the DNN, {i.e.,
\[
    \Theta = \{\bm{W}^\ell, \bm{b}^\ell, \ell=1, \cdots,L\}.
\]}
For a fully connected DNN, we have $N  =  \sum_{\ell=1}^L n_{\ell+1}(n_{\ell} + 1)$.

Our deep least squares approach uses the DNN structure~\eqref{eqn:DNN} to approximation the solution $u(\bm{x})$.  However, unlike traditional approaches~\cite{DissanayakePhan-Thien1994,CaiChenLiuLiu2019,EYu,WangZhang}, which only uses one DNN to approximate the solution $u(\bm{x})$ on the whole domain $\Omega$, we use two DNN structures to approximate $u_1(\bm{x})$ and $u_2(\bm{x})$ on $\Omega_1$ and $\Omega_2$, respectively.  In particular, for $\bm{x} \in \Omega_i$, $i=1,2$, we use DNN~\eqref{eqn:DNN} to approximation $u_i(\bm{x})$, $i=1,2$, as follows,
\begin{equation}\label{eqn:UiNN}
u_i(\bm{x}) \approx \mathcal{U}_{i, \mathcal{NN}}(\bm{x}; \Theta_i) := \bm{T}_i^L \circ \mathcal{N}_i^{L-1} \circ \cdots \circ \mathcal{N}_i^{2} \circ \mathcal{N}_i^1(\bm{x}), \quad i=1,2,
\end{equation}
with the input $\bm{x} \in \mathbb{R}^d$, i.e., $n_1 = d$ and the output $u_{i, \mathcal{NN}} \in \mathbb{R}$, i.e., $n_L = 1$. Then, the overall approximation of $u(\bm{x})$ can be defined as follows,
\begin{equation}\label{eqn:UNN}
u(\bm{x}) \approx \mathcal{U}_{\mathcal{NN}}(\bm{x}; \Theta) = 
\begin{cases}
\mathcal{U}_{1, \mathcal{NN}}(\bm{x};\Theta_1), \quad \mbox{if } \bm{x} \in \Omega_1, \\
\mathcal{U}_{2, \mathcal{NN}}(\bm{x};\Theta_2), \quad \mbox{if } \bm{x} \in \Omega_2.
\end{cases}
\end{equation}
where $\Theta = \Theta_i$ if $\bm{x} \in \Omega_i$, $i=1,2$.

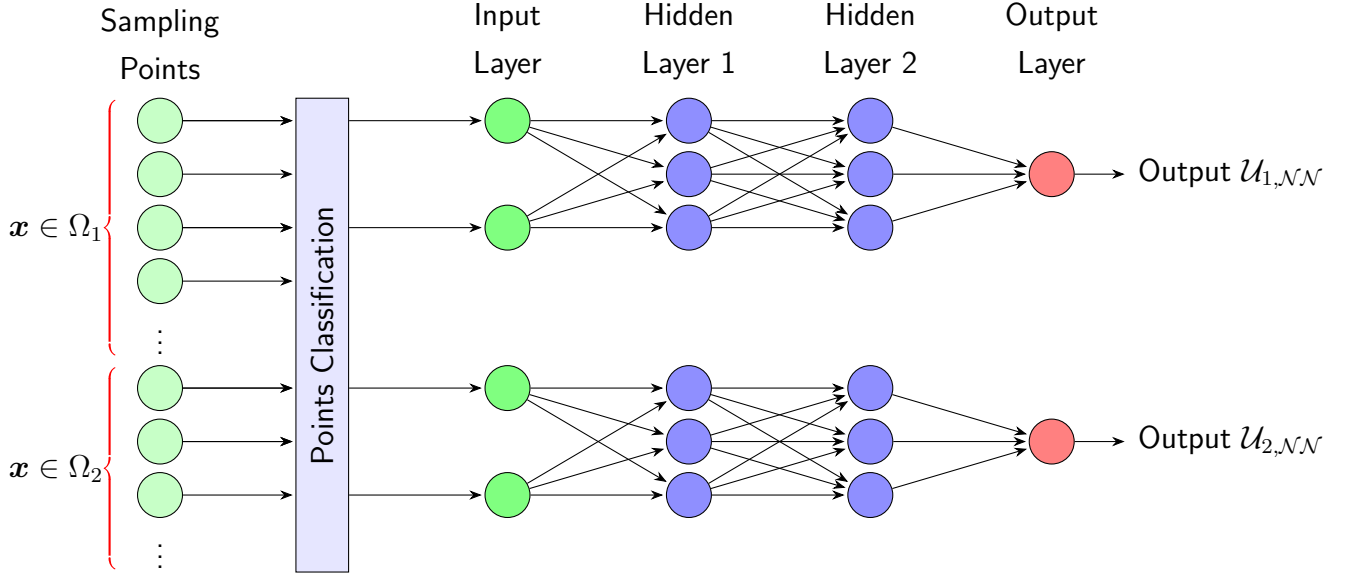
\begin{figure}[h]
\centering
\begin{tikzpicture}[
             font = \sffamily,
        shorten > = 1pt,
                > = Stealth,
    node distance = 1mm and 18mm,
      start chain = going below,
 every pin/.style = {pin distance=7mm, pin edge={thin, black, ->}},
    neuron/.style = {circle, draw, fill=#1,
                     minimum size=17pt, inner sep=.5pt,
                     },
     annot/.style = {text width=6em, align=center},
        BC/.style = {decorate,
                            decoration={calligraphic brace, amplitude=4pt,
                            pre =moveto, pre  length=1pt,
                            post=moveto, post length=1pt,
                            raise=6mm, mirror},
                            thick,
                            pen colour={red}
                         }
                        ]
\foreach \i in {1,...,9}
{
 \ifnum\i=5
    \node[neuron=white,draw=none,on chain]   (in-\i)        {$\vdots$};
 \else
     \ifnum\i<9
        \node[neuron=white,on chain,fill=green!22]             (in-\i)    {};
     \else
        \node[neuron=white,draw=none,on chain]   (in-\i)    {$\vdots$};
     \fi
 \fi
}
\path   let \p1 = ($(in-1.north) - (in-9.south)$),
            \n1 = {veclen(\y1,\x1)} in
        node (r) [minimum height=\n1,
                  minimum width=7mm, draw, fill=blue!10,
                  below right=0mm and 18 mm of in-1.north]
                  {\rotatebox{90}{Points Classification}};
\node[neuron=green!50,right=of in-1 -| r.east] (I-1) {};
\node[neuron=green!50,right=of in-3 -| r.east] (I-2) {};
\node[neuron=green!50,right=of in-6 -| r.east] (I-3) {};
\node[neuron=green!50,right=of in-8 -| r.east] (I-4) {};

\node[neuron=blue!44,right=of in-1 -| I-1.east] (H-1) {};
\node[neuron=blue!44,right=of in-2 -| I-1.east] (H-2) {};
\node[neuron=blue!44,right=of in-3 -| I-1.east] (H-3) {};

\node[neuron=blue!44,right=of in-6 -| I-1.east] (H-4) {};
\node[neuron=blue!44,right=of in-7 -| I-1.east] (H-5) {};
\node[neuron=blue!44,right=of in-8 -| I-1.east] (H-6) {};

\node[neuron=blue!44,right=of in-1 -| H-1.east] (L-1) {};
\node[neuron=blue!44,right=of in-2 -| H-2.east] (L-2) {};
\node[neuron=blue!44,right=of in-3 -| H-3.east] (L-3) {};

\node[neuron=blue!44,right=of in-6 -| H-1.east] (L-4) {};
\node[neuron=blue!44,right=of in-7 -| H-2.east] (L-5) {};
\node[neuron=blue!44,right=of in-8 -| H-3.east] (L-6) {};

\node[neuron=red!50,right=of L-2,pin=right:Output $\mathcal{U}_{1,\mathcal{NN}}$] (out) {};
\node[neuron=red!50,right=of L-5,pin=right:Output $\mathcal{U}_{2,\mathcal{NN}}$] (out2) {};
\node[annot,above=of in-1]                 {Sampling Points};
\node[annot,above=of I-1]                   {Input\\ Layer};
\node[annot,above=of I-1.north -| H-1]      {Hidden\\ Layer 1};
\node[annot,above=of I-1.north -| L-1]      {Hidden\\ Layer 2};
\node[annot,above=of I-1.north -| out]    {Output\\ Layer};

\draw[BC] (in-1.north) --
    node[left=6mm,align=right] {$\bm{x}\in\Omega_1$} (in-5.south) ;
\draw[BC] (in-6.north) --
    node[left=6mm,align=right] {$\bm{x}\in\Omega_2$} (in-9.south);
\foreach \i in {1,...,4,6,7,8}
    {\draw[->]   (in-\i) -- (in-\i -| r.west);}
\foreach \i in {1,3,6,7}
    {\draw[->]   (in-\i) -- (in-\i -| r.west);}
\foreach \i in {1,...,4}
    {\draw[->]   (I-\i -| r.east) -- (I-\i);}
\foreach \i in {1,...,2}
    \foreach \j in {1,...,3}
    \draw[->] (I-\i) edge (H-\j);
    
\foreach \i in {3,4}
    \foreach \j in {4,...,6}
    \draw[->] (I-\i) edge (H-\j);
 
\foreach \i in {1,2,3}
    \foreach \j in {1,...,3}
    \draw[->] (H-\i) edge (L-\j);
    
\foreach \i in {4,...,6}
    \foreach \j in {4,...,6}
    \draw[->] (H-\i) edge (L-\j);
    
 \foreach \i in {1,2,3}
    \draw[->] (L-\i) edge (out);
 \foreach \i in {4,5,6}
    \draw[->] (L-\i) edge (out2);
   \end{tikzpicture}
   \caption{Illustration of neural network architecture diagram}\label{fig:NN}
  \end{figure}

The DNN structure of the approximation $\mathcal{U}_{\mathcal{NN}}(\bm{x}; \Theta)$ is shown in Figure~\ref{fig:NN} for the 2 hidden layer case.  As we can see, for a sampling point, we first classify the point by its location and then determine which DNN structure to use. Such an approach gives us the freedom to approximate a solution that has severe singularities along the interface. 

With such a DNN structure, we can easily evaluate the jump along the interface as follows, for a sample point $\bm{x}_k \in \Gamma$, we have
\begin{align*}
\ljump \mathcal{U}_{\mathcal{NN}}(\bm{x}_k; \Theta) \rjump &= \mathcal{U}_{1,\mathcal{NN}}(\bm{x}_k; \Theta_1) - \mathcal{U}_{2,\mathcal{NN}}(\bm{x}_k; \Theta_2), \\
\ljump \beta(\bm{x}_k) \nabla \mathcal{U}_{\mathcal{NN}}(\bm{x}_k; \Theta) \cdot  \bm{n} \rjump &= \beta_1 \nabla \mathcal{U}_{1,\mathcal{NN}}(\bm{x}_k; \Theta_1) \cdot  \bm{n} - \beta_2 \nabla \mathcal{U}_{2,\mathcal{NN}}(\bm{x}_k; \Theta_2) \cdot  \bm{n}.
\end{align*}
These will help us to handle the interface conditions~\eqref{eqn:ic-u} and~\eqref{eqn:ic-flux} in the LS formulation when the DNN structure $\mathcal{U}_{\mathcal{NN}}(\bm{x}; \Theta)$ is used.

\subsection{Discrete Least-squares Formulations}
Next we introduce the loss function used in our deep least-sqaures approach.  Our choice is based on the LS functional defined in~\eqref{eqn:LS-obj}.  Replacing $v(\bm{x})$ with its DNN approximation $\mathcal{V}_{\mathcal{NN}}(\bm{x}; \Theta)$ defined similarly as in~\eqref{eqn:UiNN} and~\eqref{eqn:UNN}, we naturally have the following discrete LS functional,
\begin{align}
\mathcal{J}(\mathcal{V}_{\mathcal{NN}}(\bm{x};\Theta);g_j,g_f,g_D,f) &= \|-\nabla\cdot\beta(\bm{x})\nabla \mathcal{V}_{\mathcal{NN}}(\bm{x};\Theta)-f\|_{0,\Omega}^2+
{\beta_j}\|\ljump \mathcal{V}_{\mathcal{NN}}(\bm{x};\Theta)\rjump-g_j\|_{0,\Gamma}^2\notag\\
&\quad +{\beta_f}\|\ljump\beta(\bm{x})\nabla \mathcal{V}_{\mathcal{NN}}(\bm{x};\Theta)\cdot\bm{n} \rjump-g_f\|_{0,\Gamma}^2+\alpha\|\mathcal{V}_{\mathcal{NN}}(\bm{x};\Theta)-g_D\|_{0,\partial\Omega}^2. \label{eqn:discrete-LS-obj}
\end{align}
Based on the above discrete LS functional, the corresponding LS formulation is defined as follows,
\begin{equation}\label{eqn:discrete-LS-problem}
\min_{\Theta \in \mathbb{S}^N}\mathcal{J}(\mathcal{V}_{\mathcal{NN}}(\bm{x}; \Theta);g_j,g_f,g_D,f),
\end{equation}
here $\mathbb{S}^N := \{\Theta: \Theta|_{\Omega_i} \in \mathbb{R}^N, i=1,2\}$.
Standard optimization algorithms can be applied to solve the minimization problem~\eqref{eqn:discrete-LS-obj} and~\eqref{eqn:discrete-LS-problem}, which gives a straightforward way to use the DNN approximation to solve an interface problem. When only one DNN is used, this is basically the DNN methods proposed in {\cite{WangZhang}} to the interface problem~\eqref{eq:pde}-\eqref{eq:pde-bc}.

However, one difficulty in the above LS formulation is the evaluation of the norms which involves computing the integrals. In~\cite{CaiChenLiuLiu2019}, an underlying mesh is used to aid the computation of the integrals.  But, due to the existence of the interface $\Gamma$, using a mesh to compute the norms, especially the norms on the interface, is still challenging.  Therefore, we adopt a Monte-Carlo type sampling approach here and replace the discrete LS functional~\eqref{eqn:discrete-LS-obj} by the mean squared error loss function which not only helps us efficiently computing the integrals but also allows us to take advantage of advanced optimization algorithms developed in the machine learning community such as the SGD method and its variants {\cite{RobbinsMonro1951a,hardt2015train}}. 

The basic idea of our approach is to sample some points in the domain~$\Omega$ and use those sampled points to mimic the LS functional and define the discrete loss function.  In particular, we sample $M_i$ points $\{\bm{x}_k^{\Omega_i} \}_{k=1}^{M_i}\in \Omega_i$, $i=1,2$, and define the following loss function,
\begin{equation*}
\mathcal{L}_{i}(\Theta) := \frac{1}{M_i} \sum_{k=1}^{M_i}\left| -\nabla \cdot \beta_i \nabla \mathcal{U}_{i, \mathcal{NN}}(\bm{x}_k^{\Omega_i}; \Theta_i) - f(\bm{x}_k^{\Omega_i}) \right |^2,
\end{equation*}
which approximates the first term on the right-hand side of~\eqref{eqn:discrete-LS-obj}. In the practical implementation, the derivative can be replaced by backward/forward difference or performed by employing an automatic differentiation package. We also sample $M_{\partial \Omega}$ points $\{  \bm{x}_k^{\partial \Omega}  \}_{k=1}^{M_{\partial \Omega}}\subset \partial \Omega$ and approximate the boundary term (the last term) on the right hand side of~\eqref{eqn:discrete-LS-obj} as follows,
\begin{equation*}
\mathcal{L}_{\partial\Omega}(\Theta) :=
\frac{\alpha}{M_{\partial \Omega}} \sum_{k=1}^{M_{\partial \Omega}} \left| 
\mathcal{U}_{\mathcal{NN}}(\bm{x}_k^{\partial \Omega}; \Theta) - g_D(\bm{x}_k^{\partial \Omega}) \right|^2.
\end{equation*}
Finally, to handle the interface condition, we sample $M_{\Gamma}$ points $\{ \bm{x}_k^{\Gamma} \}_{k=1}^{M_{\Gamma}} \subset \Gamma$ and define the following discrete loss function on the interface,
\begin{equation*}
\mathcal{L}_{\Gamma}(\Theta) :=  \frac{{\beta_j}}{M_{\Gamma}} \sum_{k=1}^{M_{\Gamma}} \left| \ljump \mathcal{U}_{\mathcal{NN}}(\bm{x}_k^{\Gamma}; \Theta) \rjump - g_j(\bm{x}_k^{\Gamma})  \right|^2  +  \frac{{\beta_f}}{M_{\Gamma}} \sum_{k=1}^{M_{\Gamma}} \left|  \ljump \beta(\bm{x}_k^{\Gamma}) \nabla \mathcal{U}_{\mathcal{NN}}(\bm{x}_k^{\Gamma}; \Theta) \cdot \bm{n} \rjump - g_f(\bm{x}_k^{\Gamma})  \right|^2.
\end{equation*}

Now, we are ready to define the total loss function as follows,
\begin{equation}\label{eqn:total-loss}
\mathcal{L}_{\text{total}}(\Theta):= \mathcal{L}_1(\Theta) + \mathcal{L}_2(\Theta) + \mathcal{L}_{\Gamma}(\Theta) + \mathcal{L}_{\partial \Omega}(\Theta),
\end{equation}
and our deep least-squares methods  for interface problem minimize the above discrete loss function~\eqref{eqn:total-loss} as follows,
\begin{equation}\label{eqn:deep-LS}
\min_{\Theta \in \mathbb{S}^N} \mathcal{L}_{\text{total}}(\Theta).
\end{equation}
Let $\Theta^*$ denote the minimizer and the corresponding DNN approximation is given by $\mathcal{U}_{\mathcal{NN}}(\bm{x}; \Theta^*)$.

\begin{remark}
In the definition of the loss function, we could weight each loss functions differently.  However, for the sake of simplicity, we use the fixed weights here, i.e., $\beta_j = \beta_f =1$ and $\alpha = 500$, and the numerical experiments show that this choice works well in practice.
\end{remark}


\section{Numerical Examples}\label{sec:numerics}

\begin{figure}[H]
	\centering
	\begin{tabular}{c}
		\includegraphics[width=.45\textwidth]{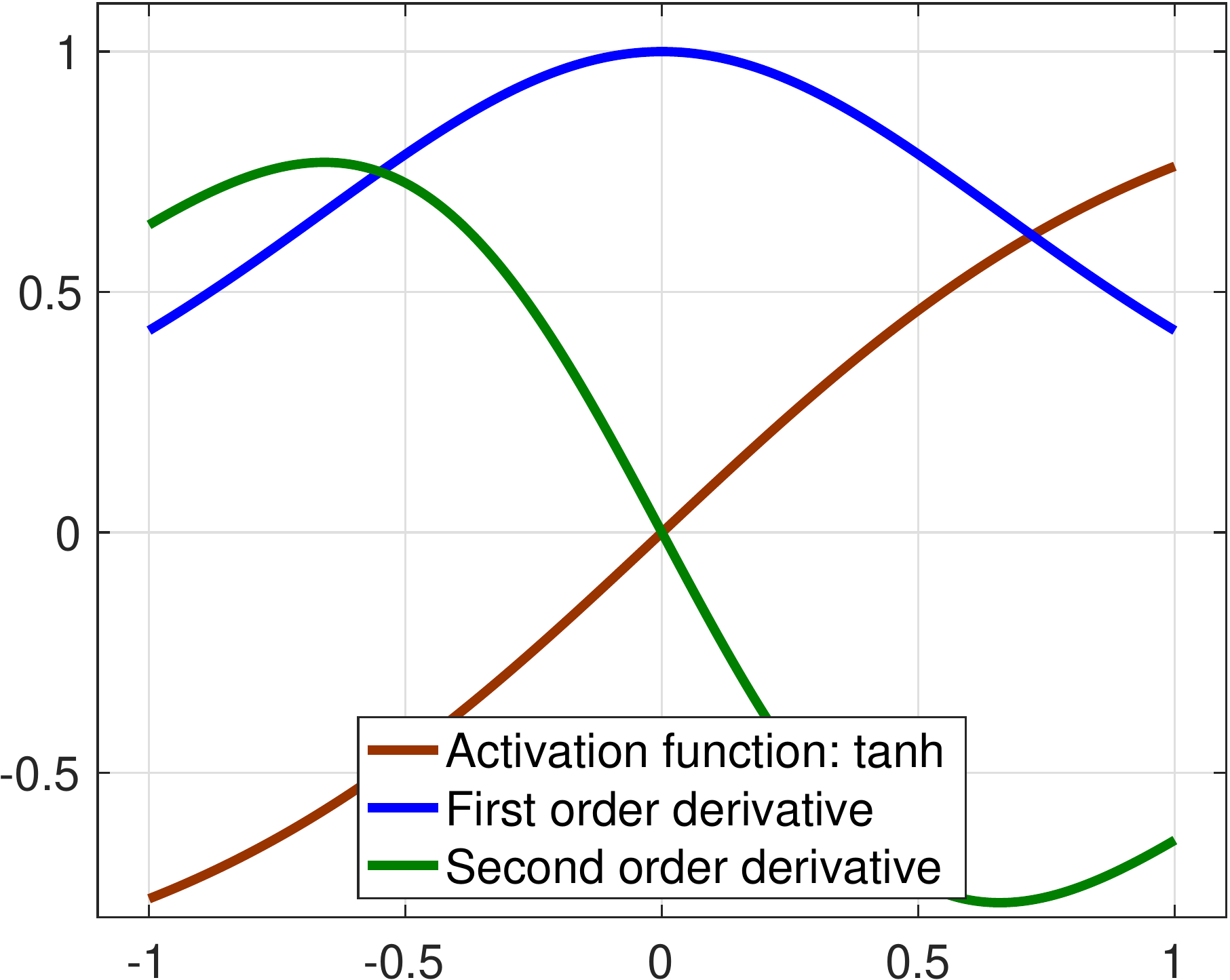}
	\end{tabular}
	\caption{Illustration of the activation function $\psi = \tanh$. }\label{Fig:tanh}
\end{figure}

In this section, 
we apply our algorithm using piecewise DNN structure based on the loss function defined in \eqref{eqn:total-loss}
to solve the elliptic interface problem~\eqref{eq:pde}-\eqref{eq:pde-bc}. 
Our numerical experiments are implemented based on TensorFlow {\cite{abadi2016tensorflow}}.  In all the examples, we choose the activation function to be $\psi = \text{tanh}$ (see Figure~\ref{Fig:tanh}). Recall the convention that $\Omega_1$ is the interior sub-domain and $\Omega_2$ is the exterior sub-domain, we shall set up two DNN structures $\mathcal{U}_{1,\mathcal{NN}}$ and $\mathcal{U}_{2, \mathcal{NN}}$ to approximate the exact solution $u_1$ and $u_2$ in $\Omega_1$ and $\Omega_2$, respectively. These two DNN structures are independent and can be set up and trained separately if needed. For the numerical experiments in this paper, we use the same number of layered neural network for both $\mathcal{U}_{1,\mathcal{NN}}$ and $\mathcal{U}_{2, \mathcal{NN}}$. In each DNN structure, a fully connected DNN is implemented. All parameters of the DNNs are trained simultaneously with the single discrete loss function defined in~\eqref{eqn:total-loss}. We choose $\alpha = 500$ in all the experiments and each layer of the DNN contains $64$ neurons. In the training process, a variant of the stochastic gradient descent method, ADAM~\cite{KingmaAdam}, is applied with an initial learning rate of $0.001$ and $2\times 10^5$ epochs. 




\subsection{Example 1. Sunflower Shape Interface}
\label{Sect:NumTest4}
In this example, we consider a sunflower-shaped interface $\Gamma$ that has parametric form as follows,
\begin{eqnarray*}
\begin{cases}
x(t) = r(\theta)\cos(\theta)+x_c,\\
y(t) = r(\theta)\sin(\theta)+y_c,
\end{cases}
\end{eqnarray*}
where $r(\theta) = r_0+r_1\sin(\omega\theta)$, $0\le\theta< 2\pi.$ The level set function is described as follows:
\begin{eqnarray*}
(x-x_c)^2+(y-y_c)^2 = r(\theta)^2. 
\end{eqnarray*}
We choose $r_0 = 0.4,$ $r_1 = 0.2$, $\omega = 20$, and $x_c = y_c = 0.02\sqrt{5}$ in our experiments. The coefficient $\beta$ is a piece-wise constant with $\beta_1 = 1$ in $\Omega_1$ and $\beta_2=10$ in $\Omega_2$. The exact solution is chosen as
\begin{eqnarray*}
u(\bm{x}) = \begin{cases}
\dfrac{r^2}{\beta_1}, &\mbox{ if } \bm{x}\in\Omega_1,\\[7pt]
\dfrac{r^4-0.1\ln(2r)}{\beta_2}, &\mbox{ if } \bm{x}\in\Omega_2.
\end{cases}
\end{eqnarray*}
The jump conditions $g_j$ and $g_f$ are then computed by the exact solution and $\beta$.

\begin{figure}[H]
\centering
\begin{tabular}{ccc}
\includegraphics[width=.32\textwidth]{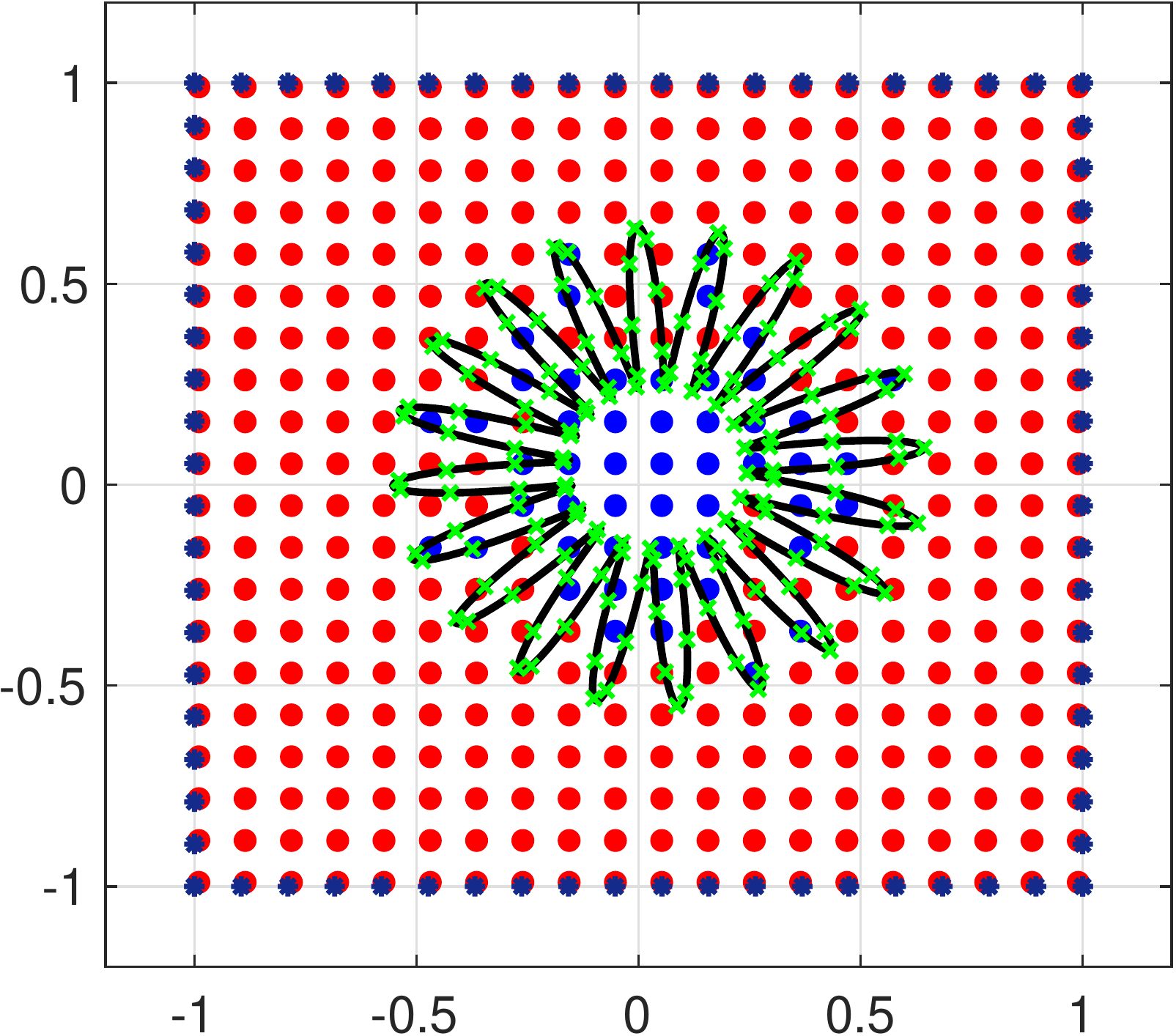}
&\includegraphics[width=.32\textwidth]{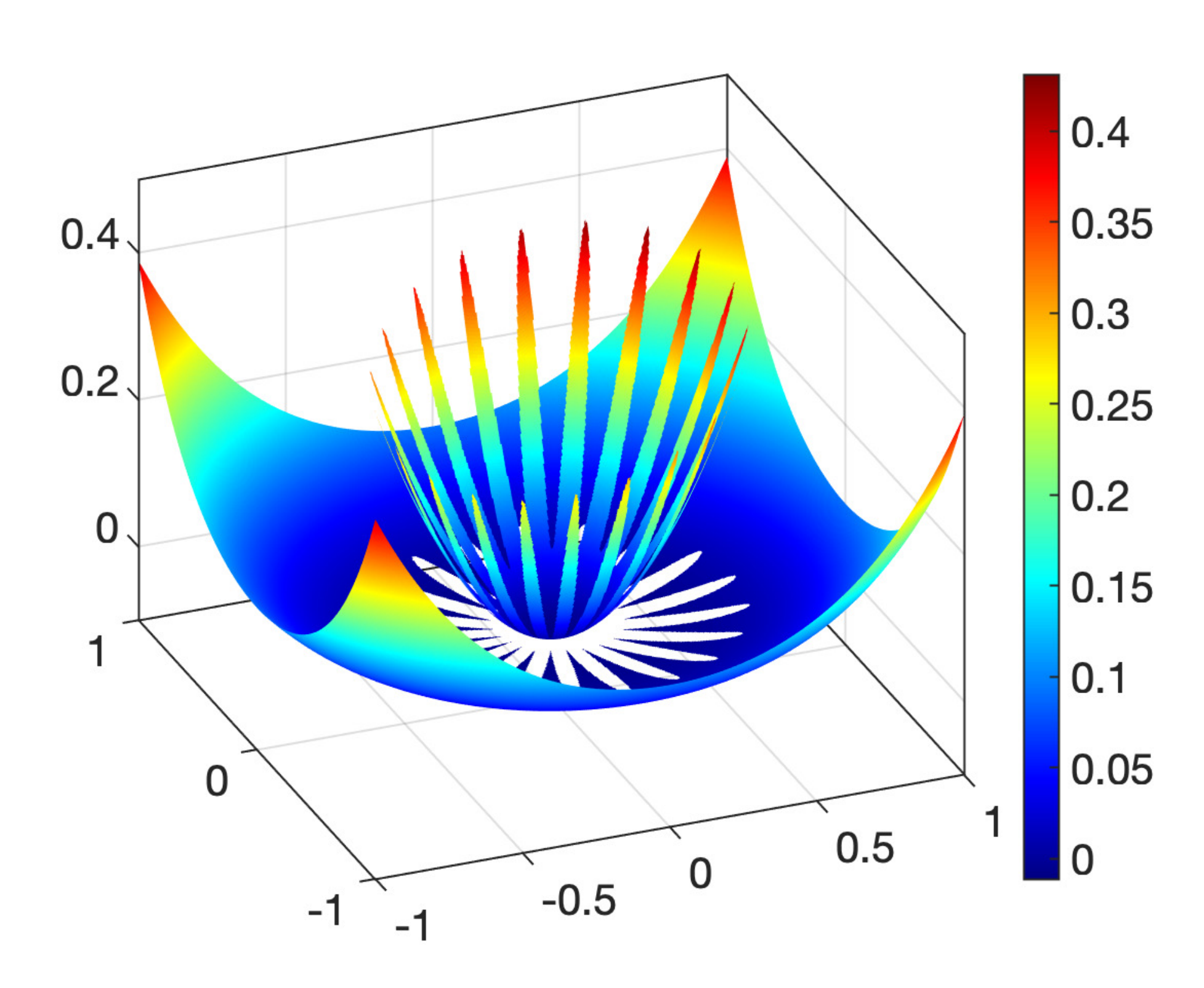}
&\includegraphics[width=.32\textwidth]{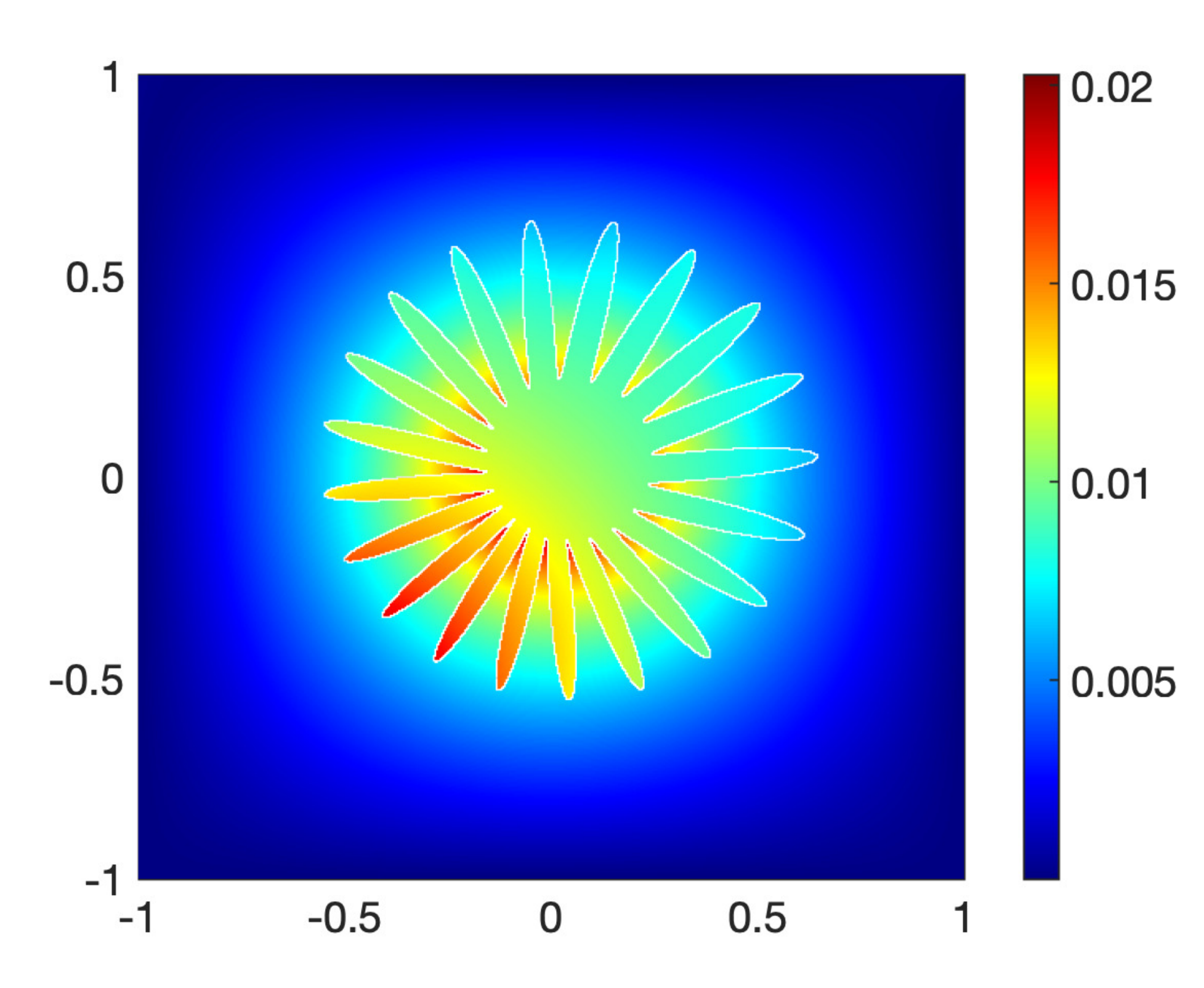}\\
(a) & (b) &(c)
\end{tabular}
\caption{ Example~\ref{Sect:NumTest4}. (a) $20\times 20$ grid;
(b). Numerical solution $\mathcal{U_{NN}}$; 
(c). Error of $u - \mathcal{U_{NN}}$
 }\label{Fig:Test4-1}
\end{figure}

We fix the neural network with $8$ layers. The numerical solution is calculated on the uniform sampled points with $M_1 = 51$ on the domain $\Omega_1$, $M_2 = 349$ on the domain $\Omega_2$ and $M_{\Gamma} = 160$ on the interface $\Gamma$, and  $M_{\partial \Omega}=80$ on the boundary $\partial \Omega$, see Figure~\ref{Fig:Test4-1}a.

The neural network approximation~$\mathcal{U}_{\mathcal{NN}}$ is plotted in Figure~\ref{Fig:Test4-1}(b) and the error $u - \mathcal{U}_{\mathcal{NN}}$, is shown in Figure~\ref{Fig:Test4-1}(c). The relative error for this case is $\dfrac{\|u - \mathcal{U}_{\mathcal{NN}}\|_{0,\Omega}}{\|u\|_{0,\Omega}}=5.3183\text{E}-2$. This example shows that our neural network algorithm is able to provide satisfactory approximations for interface problems with complex interface even on uniformly sampled points.

\subsection{Example 2. Sphere Shape Interface} \label{Sect:Num-Test5}
\begin{figure}[H]
\centering
\begin{tabular}{c}
\includegraphics[width=.45\textwidth,height = .4\textwidth]{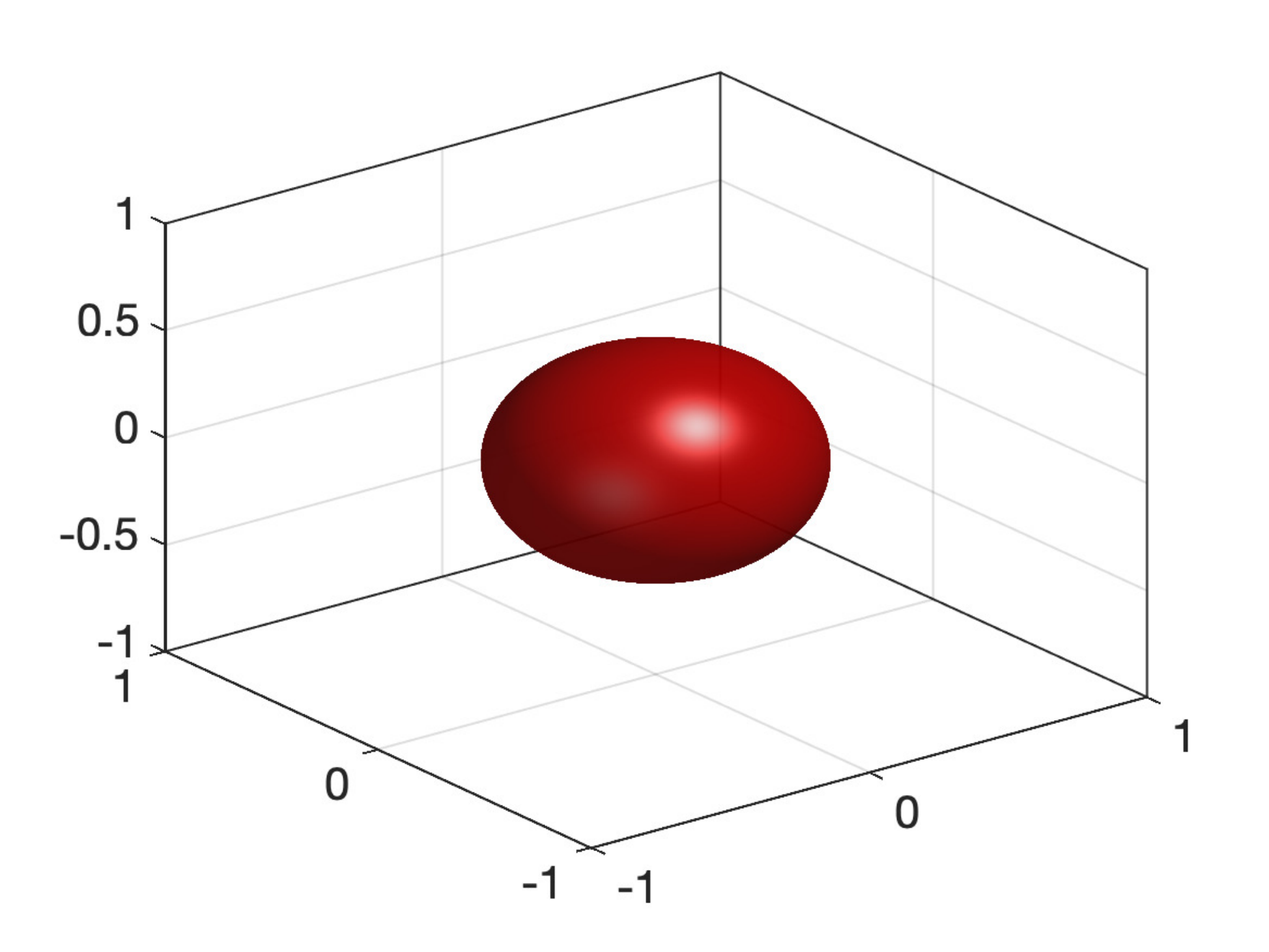}
\end{tabular}
\caption{Example~\ref{Sect:Num-Test5}. Sphere-shaped interface.
}\label{Fig:Test5-Test6-interface}
\end{figure}

\begin{figure}[H]
\centering
\begin{tabular}{ccc}
\includegraphics[width=.3\textwidth]{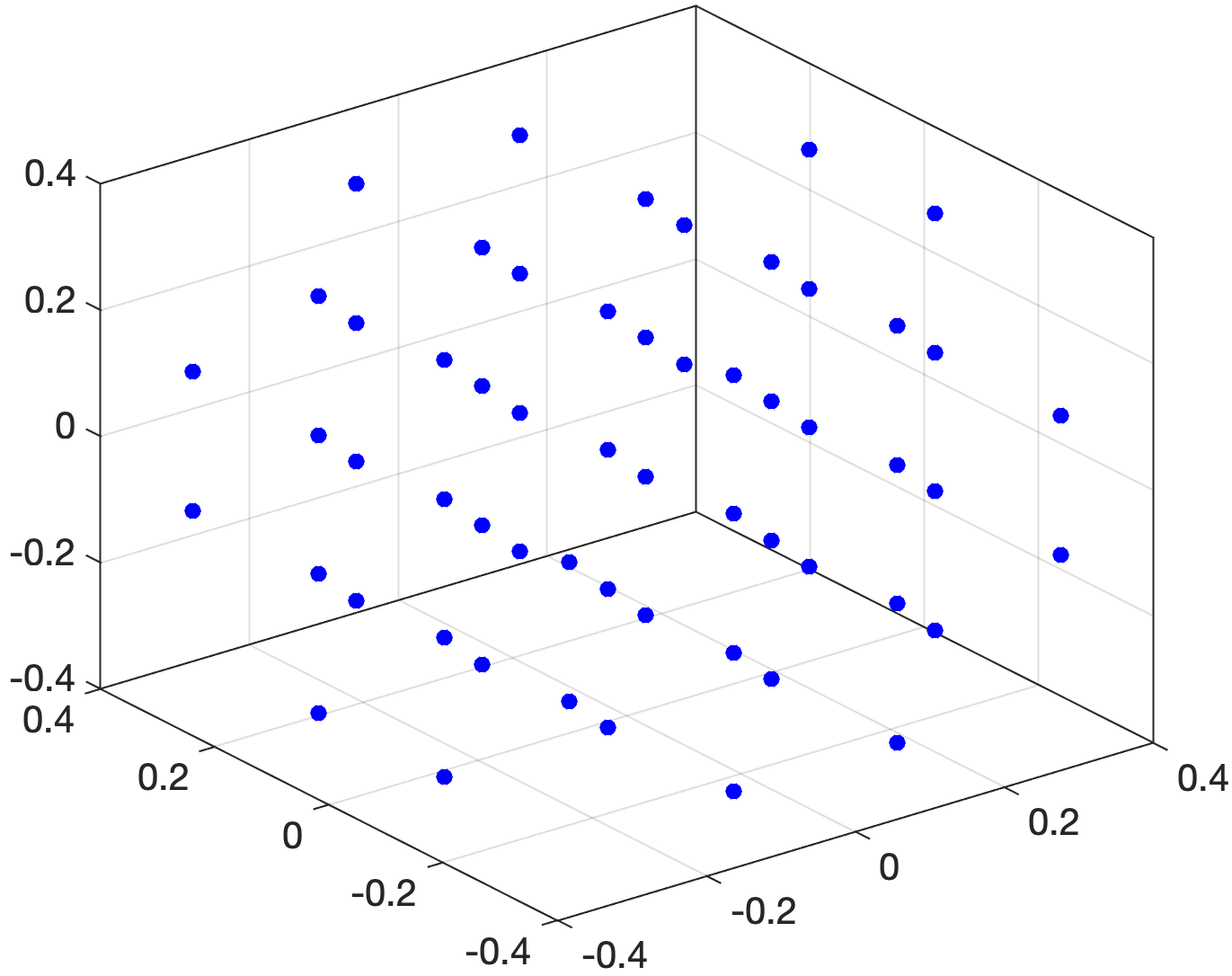}
&
\includegraphics[width=.3\textwidth]{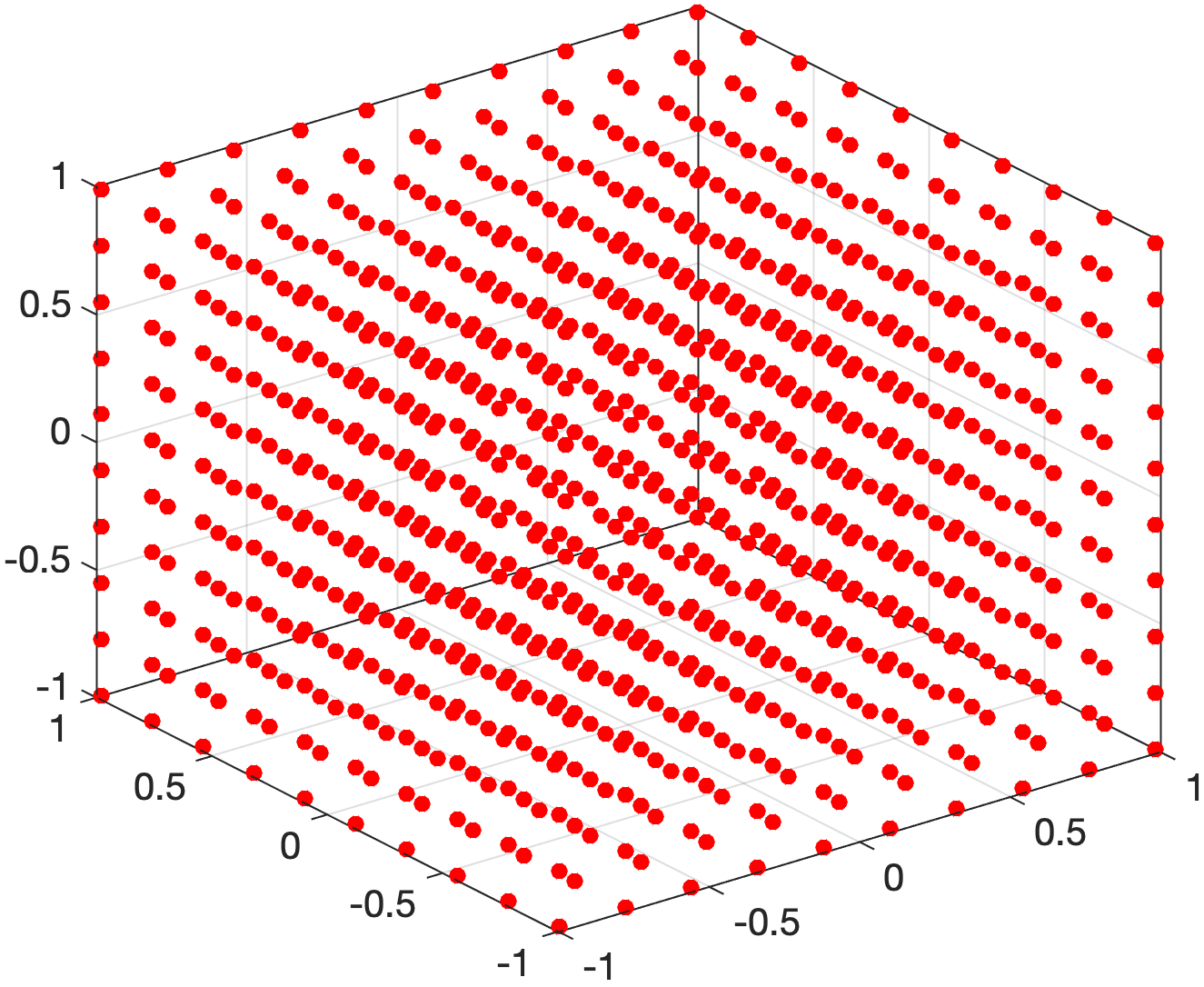}
&
\includegraphics[width=.3\textwidth]{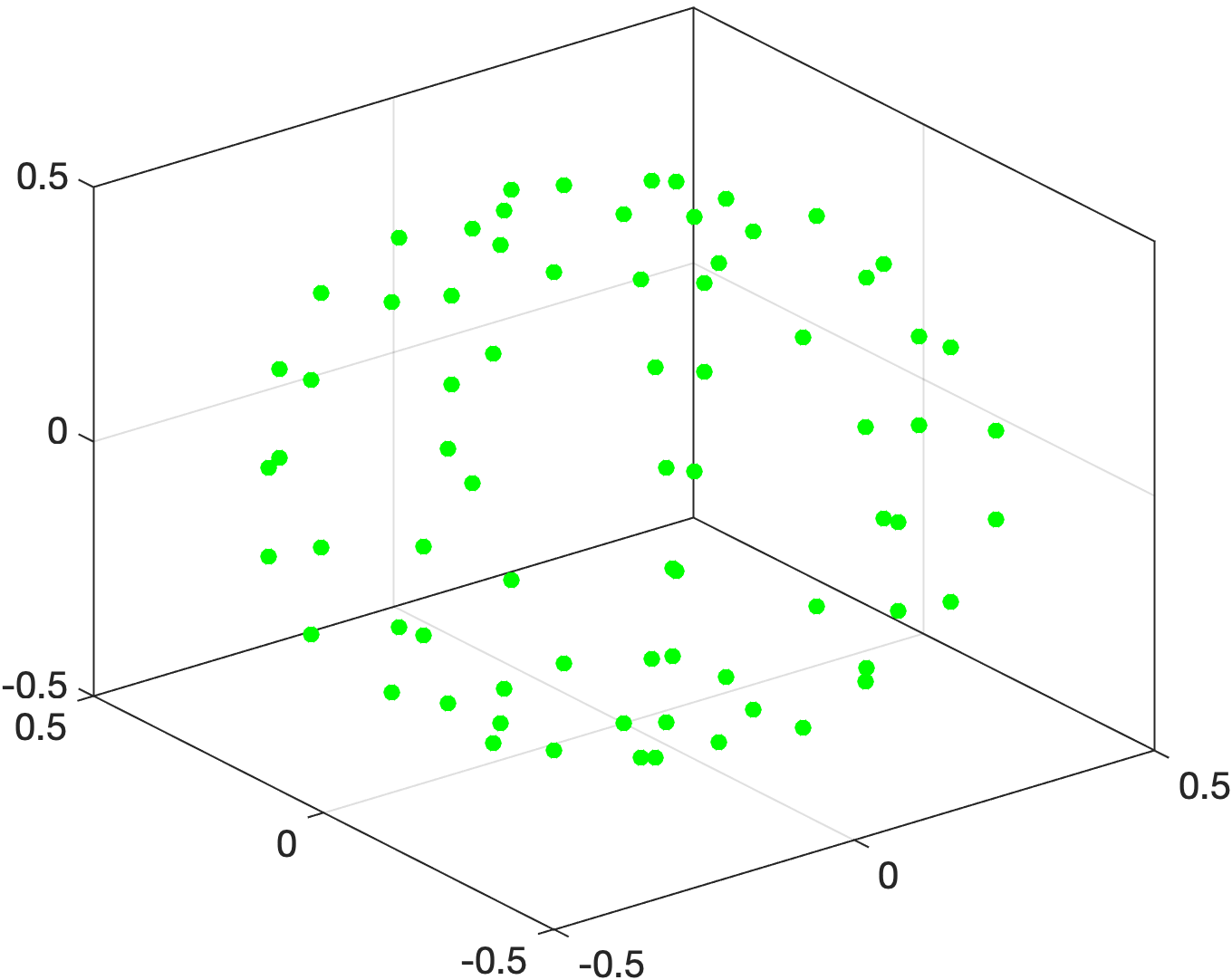}
\\
(a) & (b) & (c)
\end{tabular}
{
\caption{ Example~\ref{Sect:Num-Test5}. Sampling points in (a) $\Omega_1$; (b) $\Omega_2$; (c) $\Gamma$.}}
\label{Fig:Test5-points}
\end{figure}

In this example, we test our algorithm for a three-dimensional problem. Let $\Omega=[-1,1]^3$ and the interface is defined as the zero level set of the following level set function,
$$
\phi(x,y,z)=0.5^2 - x^2+y^2+z^2.
$$ 
The solution in two different subdomains is chosen as
\begin{eqnarray*}
u_1(x,y,z) = \cos(x)\cos(y)\cos(z)\mbox{ and }u_2 (x,y,z) = 0.
\end{eqnarray*}
The discontinuous coefficients are given by
$\beta_1 = 10$ in $\Omega_1$ and $\beta_2 = 1$ in $\Omega_2$.
We test our algorithm with $4$ hidden layers. The number of uniformly sampled points are $M_1 = 56$, $M_2 = 944$, 
$M_{\Gamma} = 100$ on the interface, and $80$ points on each face of the boundary, i.e. $M_{\partial \Omega} =480$, as shown in Figure~\ref{Fig:Test5-points}.

\begin{figure}[H]
\centering
\begin{tabular}{cc}
\includegraphics[width=.45\textwidth]{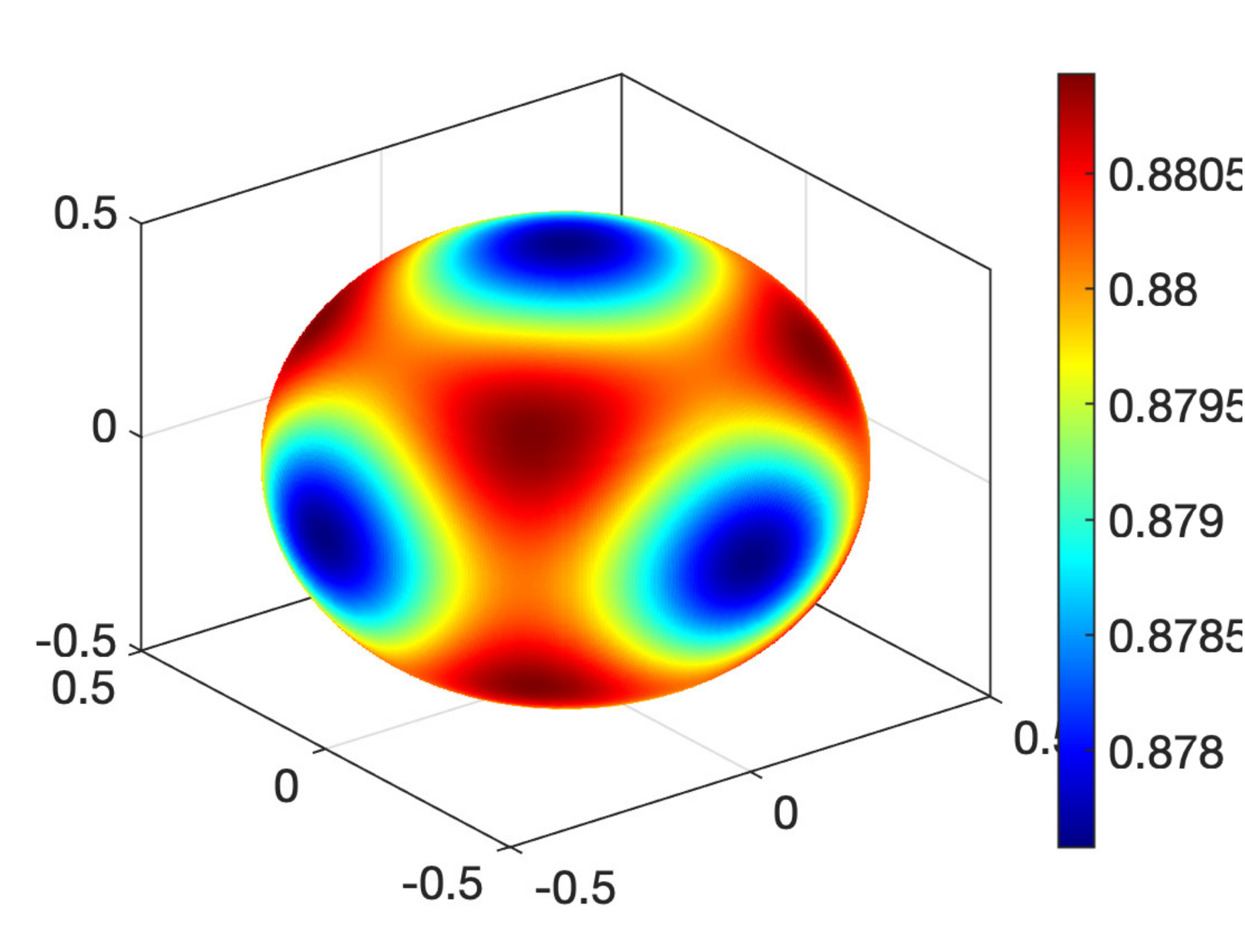}
&\includegraphics[width=.45\textwidth]{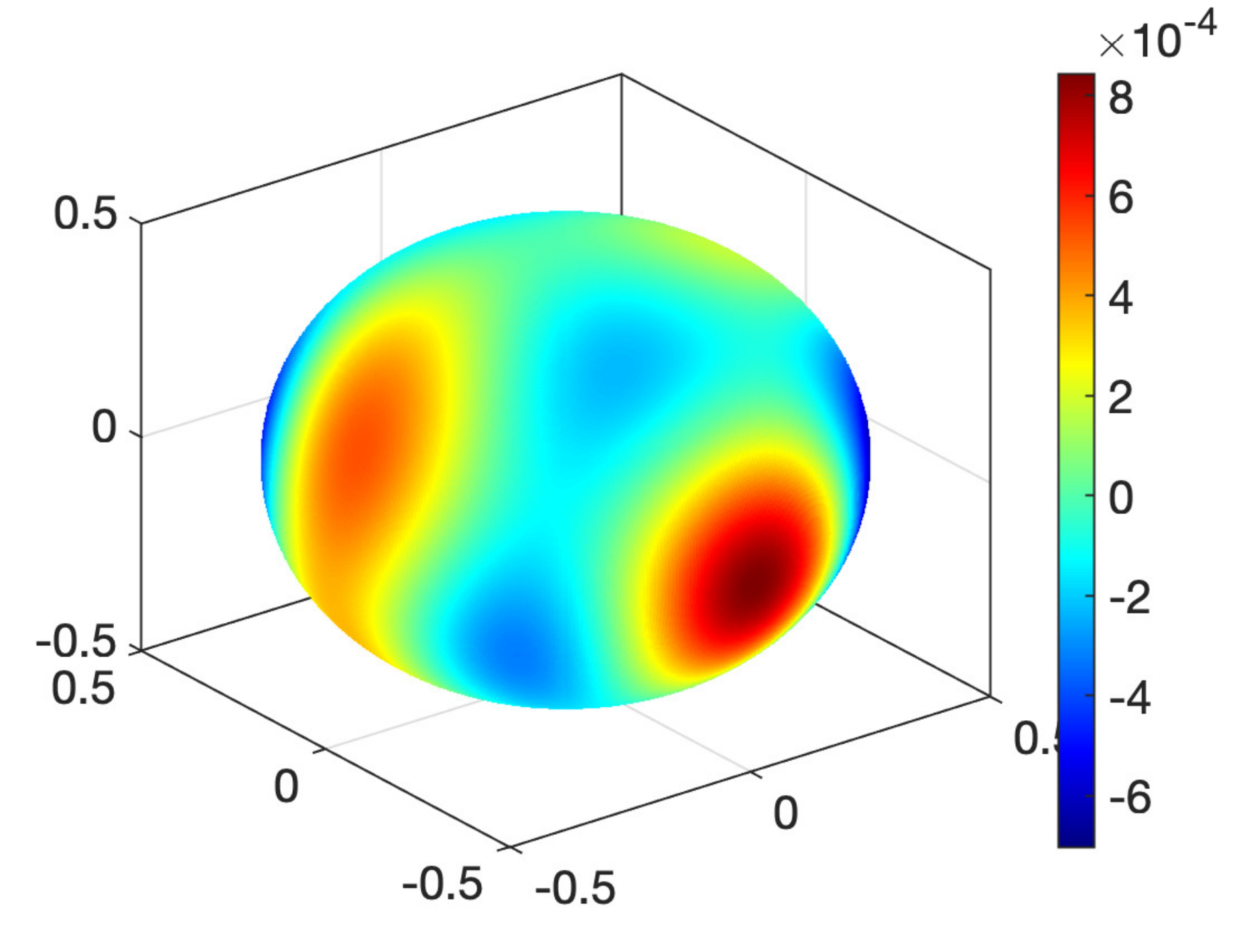}\\
(a) & (b)
\end{tabular}
\caption{Example~\ref{Sect:Num-Test5}. (a) 
Plot of $\mathcal{U_{NN}}$ on the interface $\Gamma$; (b) Plot of error, 
$u - \mathcal{U_{NN}}$, on the interface $\Gamma$. }\label{Fig:Test5-1}
\end{figure}

The numerical solution and the error,$u - \mathcal{U_{NN}}$, on the interface are plotted in Figure~\ref{Fig:Test5-1}(a) and Figure~\ref{Fig:Test5-1}(b), respectively.
The relative error in the $L^2$ norm is $\dfrac{\|u-\mathcal{U_{NN}}\|_{0,\Omega}}{\|u\|_{0,\Omega}} = 5.4508\text{E}$-5.
This example shows that our algorithm also works reasonably well for three-dimensional problems on the uniformly sampled points.

\subsection{Example 3. Heart Shape Interface}\label{Sect:Num-Test7}
\begin{figure}[H]
\centering
\includegraphics[width=.49\textwidth]{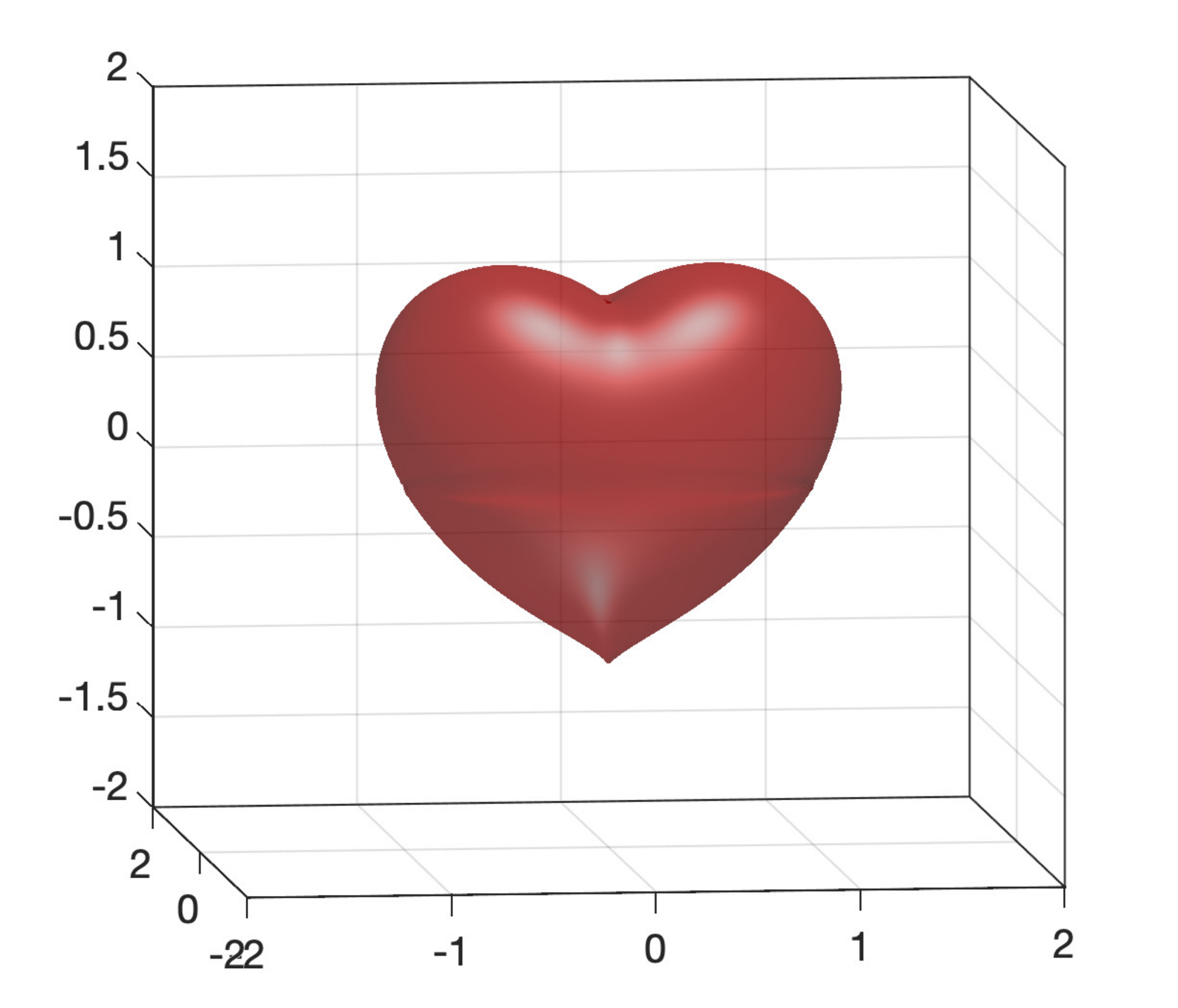}
\caption{ Example~\ref{Sect:Num-Test7}. The heart interface~\eqref{eq:Test7-interface}. }\label{Fig:Test7-interface}
\end{figure}

\begin{figure}[H]
\centering
\begin{tabular}{ccc}
\includegraphics[width=.3\textwidth]{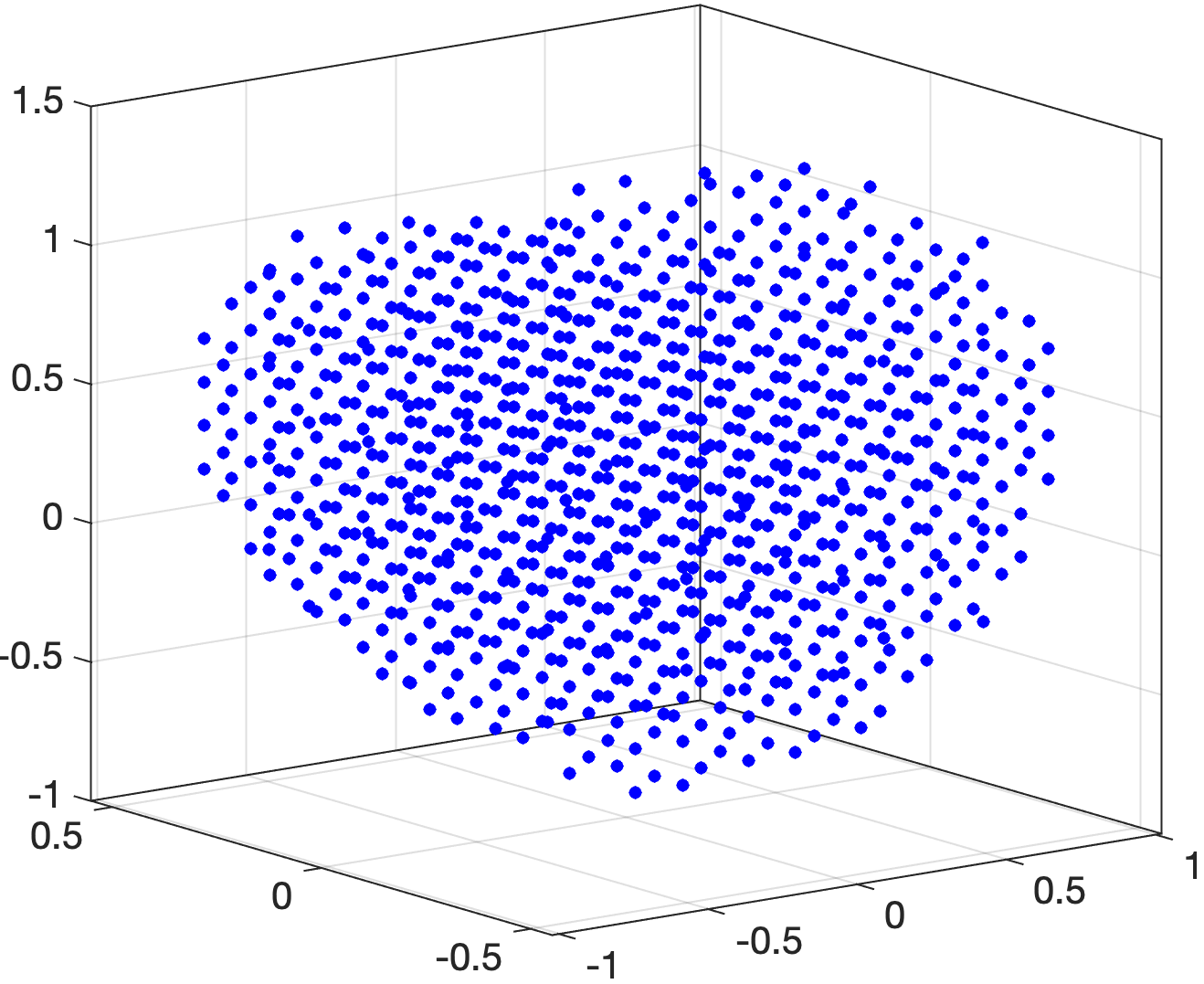}
&
\includegraphics[width=.3\textwidth]{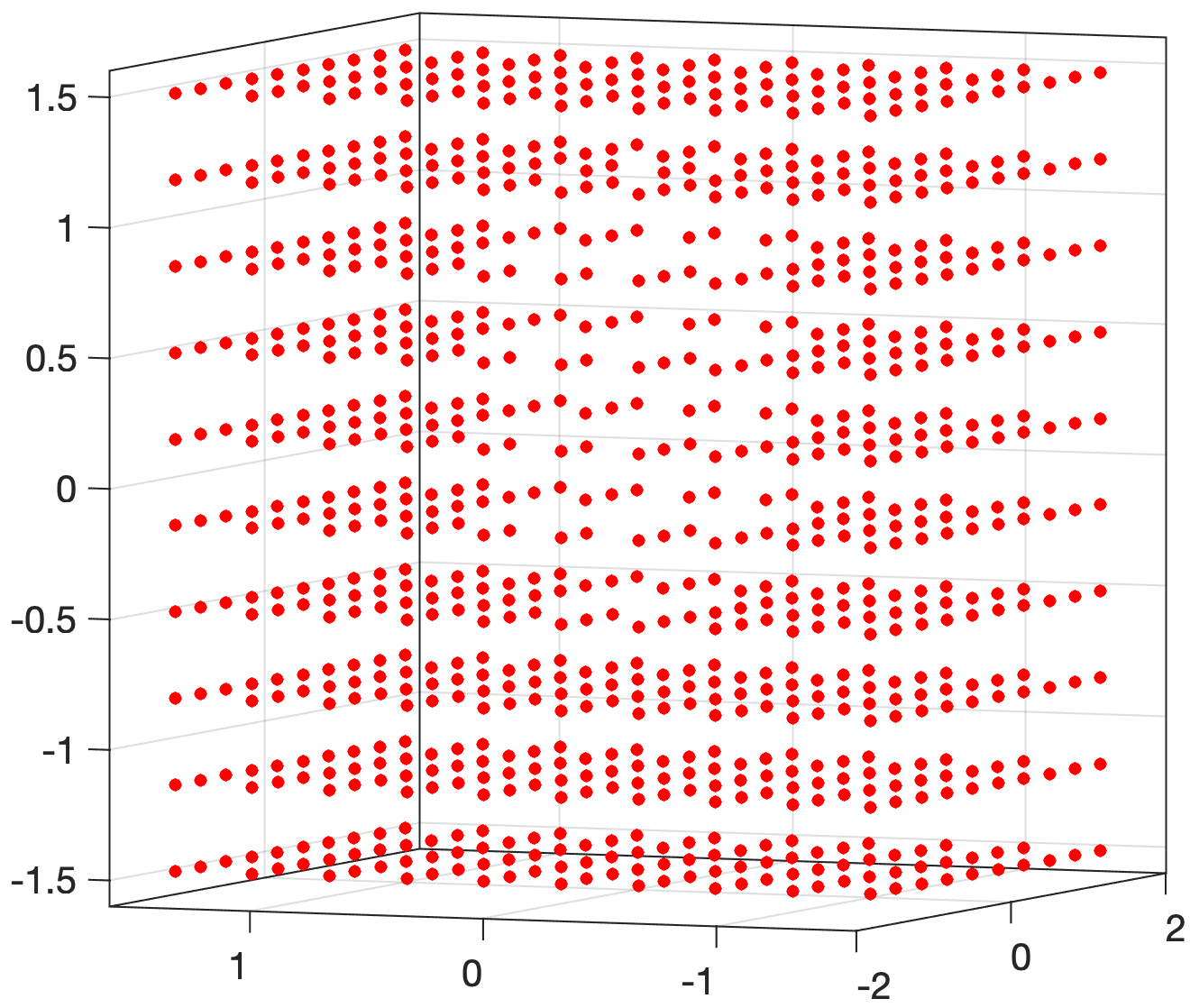}
&
\includegraphics[width=.3\textwidth]{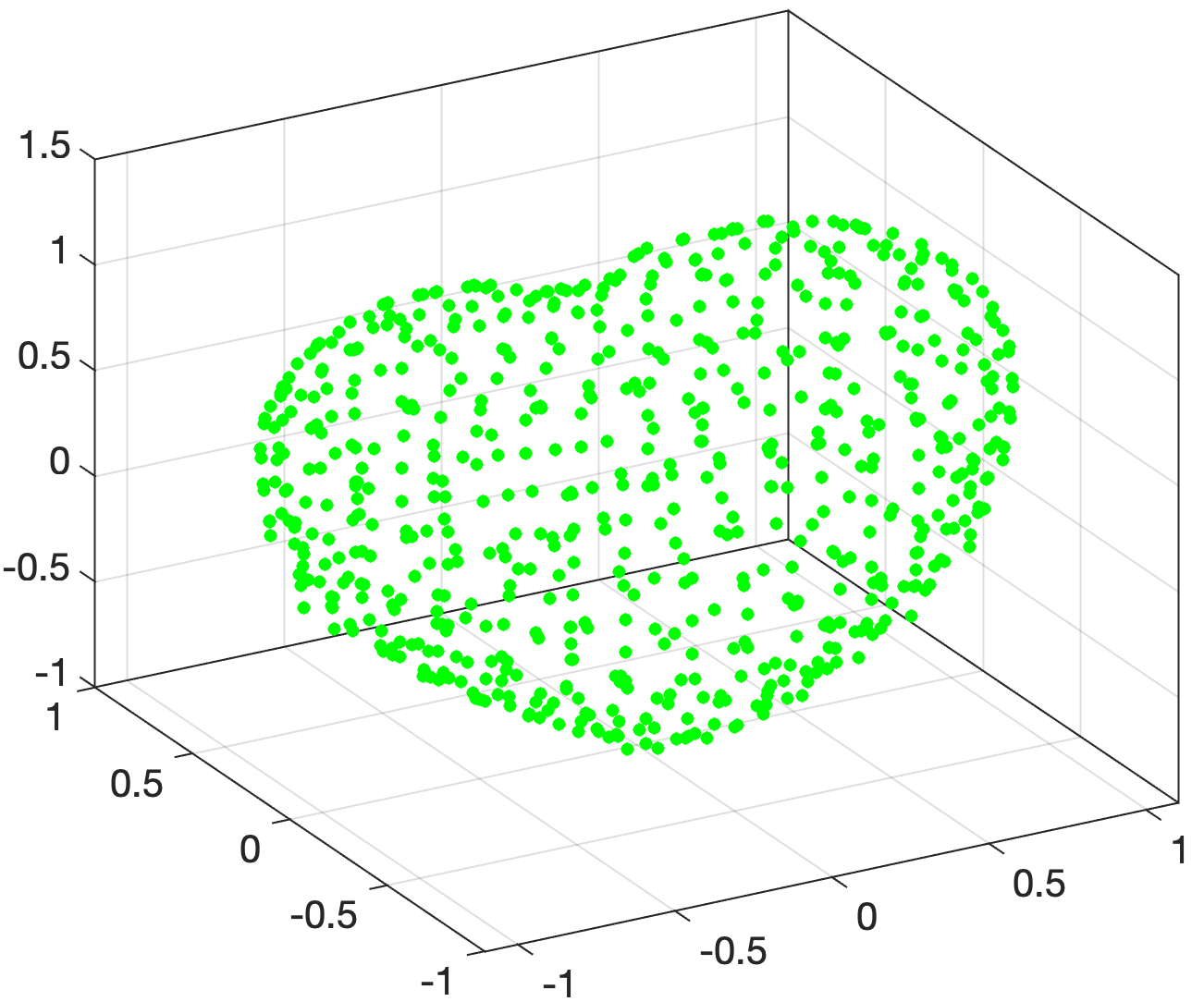}
\\
(a) & (b) & (c)
\end{tabular}
\caption{ Example~\ref{Sect:Num-Test7}. 
Sampling points in (a) $\Omega_1$; (b) $\Omega_2$; (c) $\Gamma$.}\label{Fig:Test7-points}
\end{figure}

\begin{figure}[H]
\centering
\begin{tabular}{cc}
\includegraphics[width=.46\textwidth]{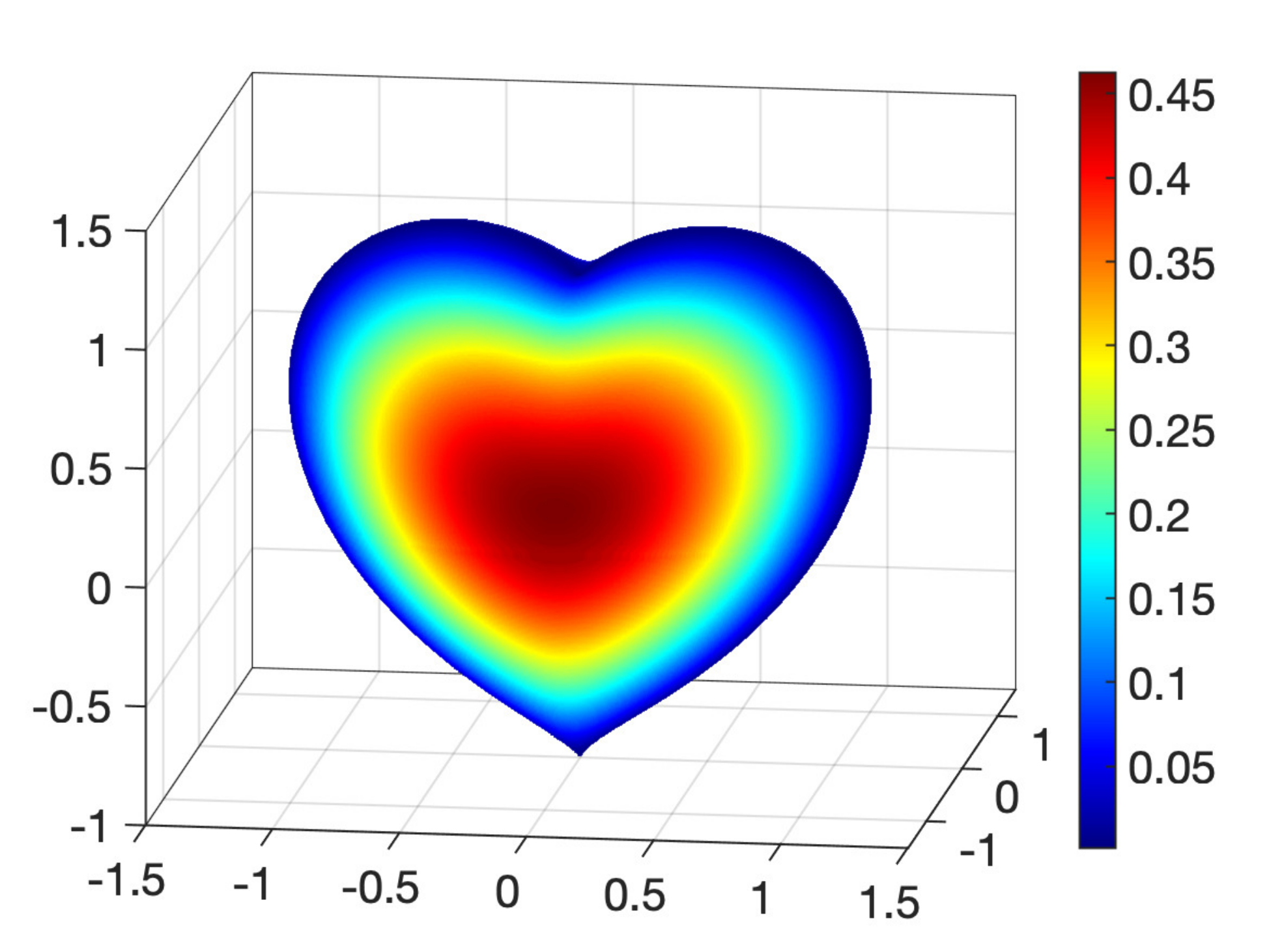}
&
\includegraphics[width=.46\textwidth]{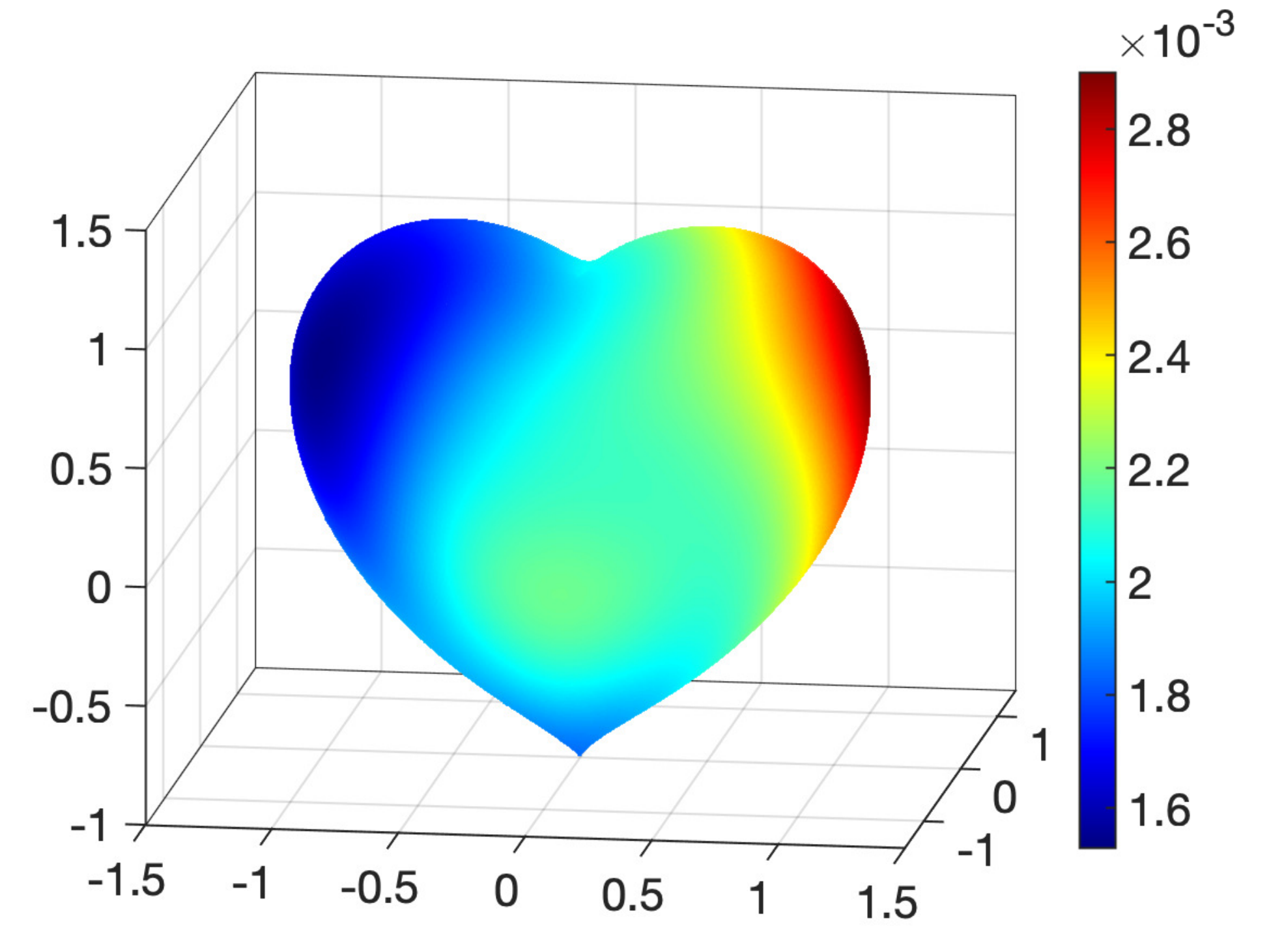}\\
(a) & (b)
\end{tabular}
\caption{ Example~\ref{Sect:Num-Test7}. (a) Plot of $\mathcal{U_{NN}}$ on the interface $\Gamma$; (b) Plot error of $u - \mathcal{U_{NN}}$ on the interface $\Gamma$.}\label{Fig:Test7-1}
\end{figure}

In this test, we test a more complicated three-dimensional interface problem with a heart shape interface, see Figure \ref{Fig:Test7-interface}. Let the domain $\Omega=[-1.5,1.5]^3$ and the interface is described as the zero level set of the following level set function,
\begin{eqnarray}
\phi(x,y,z) = \left(x^2+\frac{9}{4}y^2+z^2-1\right)^3-x^2z^3-\frac{9}{80}y^2z^3.
\label{eq:Test7-interface}
\end{eqnarray}
The exact solutions are chosen as
\begin{eqnarray}
u_1 = y^2,\text{ and }u_2 = \cos(x)\cos(y)\cos(z),
\end{eqnarray}
and the diffusion coefficients are chosen as $\beta_1 = 8$ in $\Omega_1$ and $\beta_2 = 1$ in $\Omega_2$. We test our algorithm with 8 hidden layers. The number of uniformly sampled points are $M_1 = 956$, $M_2 = 908$, 
$M_{\Gamma} = 676$ on the interface, and $80$ points on each face of the boundary (again $M_{\partial \Omega} = 480$), as shown in Figure~\ref{Fig:Test7-points}. The approximation $\mathcal{U_{NN}}$ and the error, $u - \mathcal{U_{NN}}$, on the interface are plotted in Figure~\ref{Fig:Test7-1}. The relative error in the $L^2$ norm is $\dfrac{\|u - \mathcal{U_{NN}}\|_{0,\Omega}}{\|u\|_{0,\Omega}} = 1.1520\text{E}-2$.  This shows the capability of our algorithm to tackle three-dimensional interface problems with complex interfaces reasonably well on uniformly sampled points.


\subsection{Example 5. Circle Interface with High Contrast Coefficients}\label{Sect:NumTest1}
In this example, we consider the interface problem with high contrast diffusion coefficients in (\ref{eq:pde})-(\ref{eq:pde-bc}) with homogeneous jump conditions. The exact solution is
\begin{eqnarray}
u(\bm{x})=\begin{cases}
u_1(\bm{x}) = \dfrac{r^3}{\beta_1},\mbox{ if } \bm{x} \in\Omega_1\\[7pt]
u_2(\bm{x}) = \dfrac{r^3}{\beta_2}+(\dfrac{1}{\beta_1}-\dfrac{1}{\beta_2})r_0^3,\mbox{ if }\bm{x}\in\Omega_2,
\end{cases}
\end{eqnarray}
where $\Omega_1 = \{\bm{x}\,|\,|\bm{x}|<0.5\}$, $\Omega_2 = \Omega \backslash \Omega_1$, $\Omega = [-1,1] \times [-1,1]$, and $r=\sqrt{x_1^2+x_2^2}$.
The exact interface is the zero level set of the following level set function $$\phi(\bm{x}) = x_1^2+x_2^2-(0.5)^2.$$

\begin{figure}[H]
	\centering
	\begin{tabular}{cc}
		\includegraphics[width=.45\textwidth]{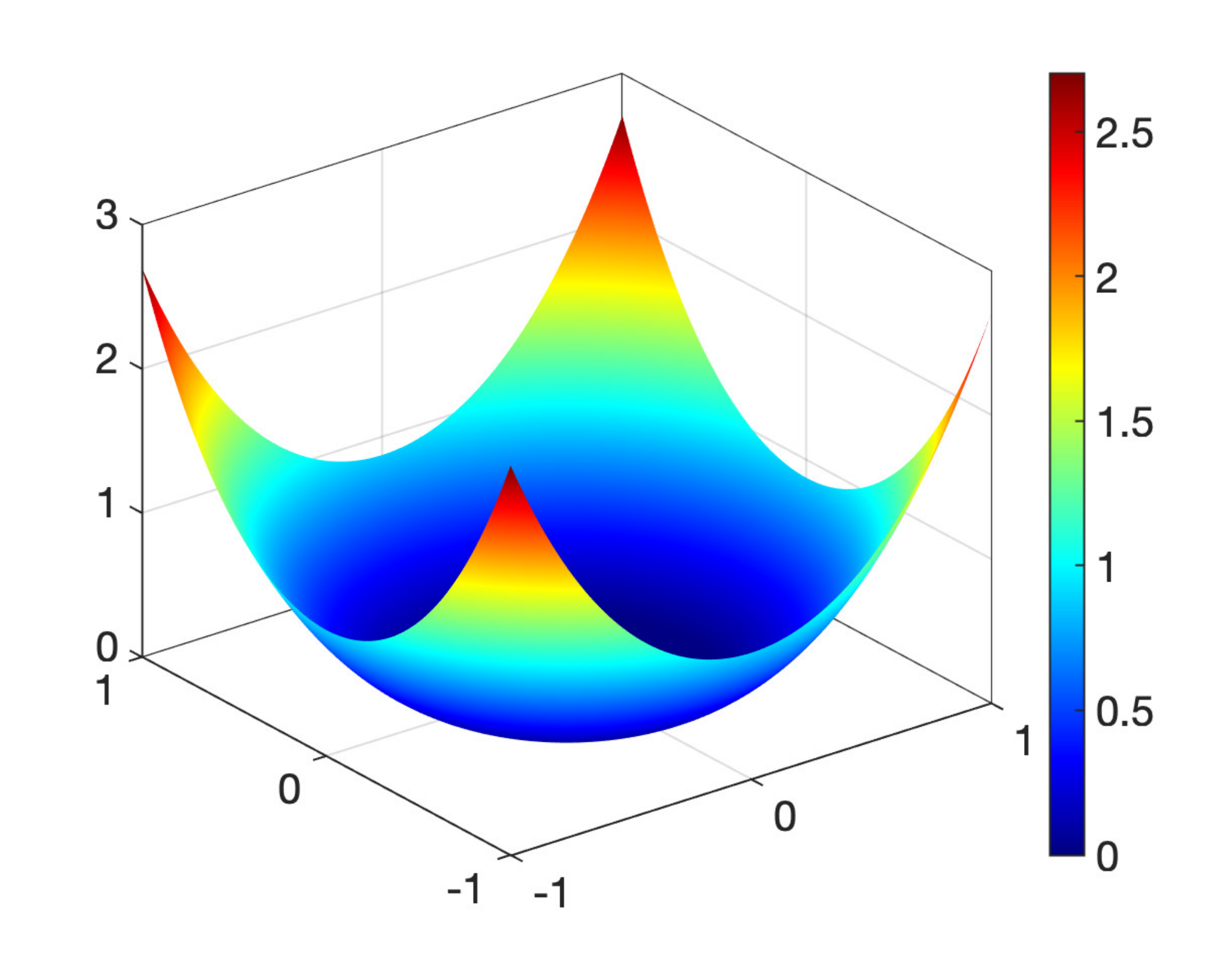}
		&\includegraphics[width=.45\textwidth]{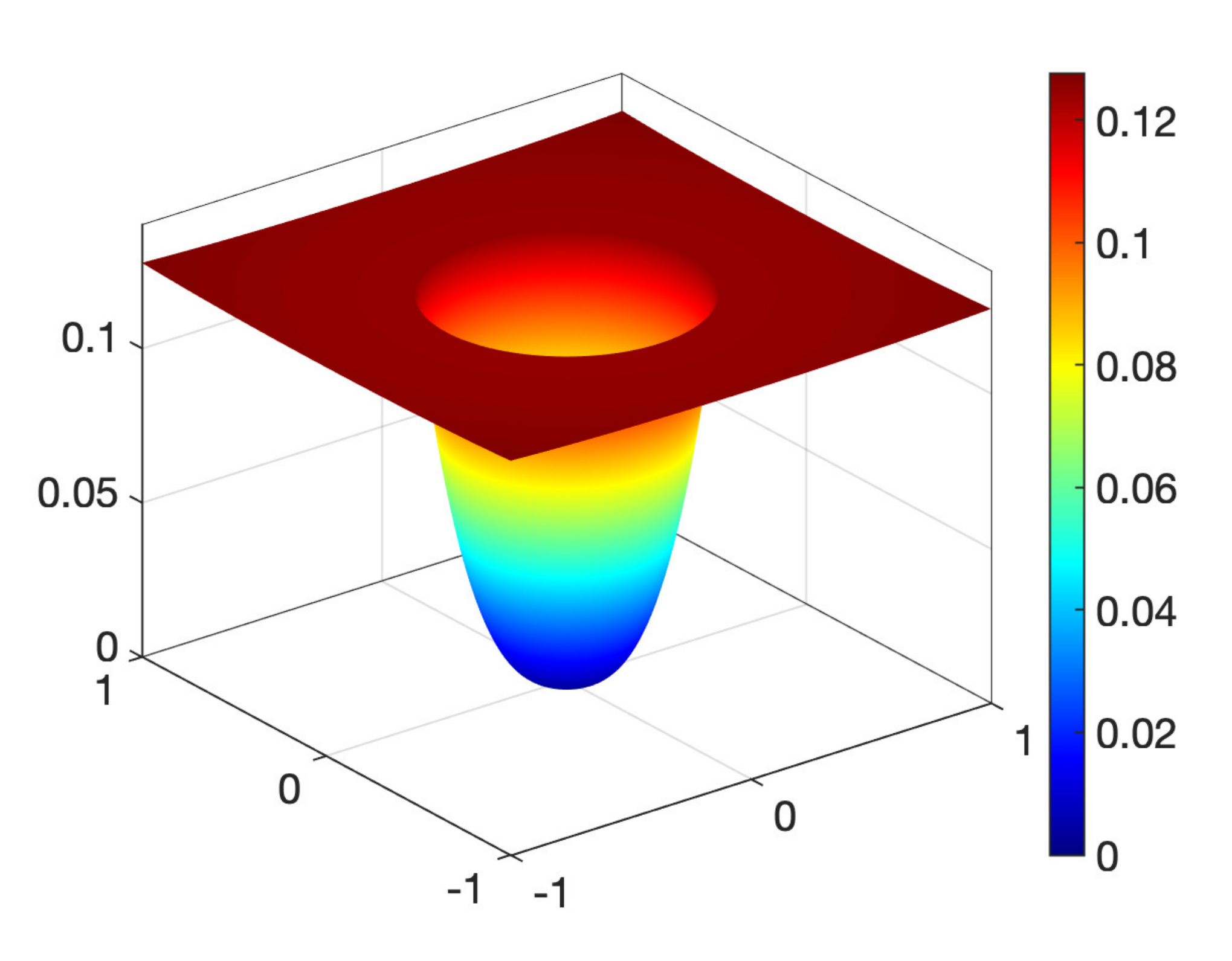}\\
		(a) & (b)
	\end{tabular}
	\caption{ Example~\ref{Sect:NumTest1}. Plots of exact solutions for (a) $\beta_1 = 1000$ and $\beta_2 = 1$; (b) $\beta_1 = 1$ and $\beta_2 = 1000$.
	}
	\label{Fig:Test1-1}
\end{figure}

\begin{figure}[H]
	\centering
	\begin{tabular}{ccc}
		\includegraphics[width=.32\textwidth]{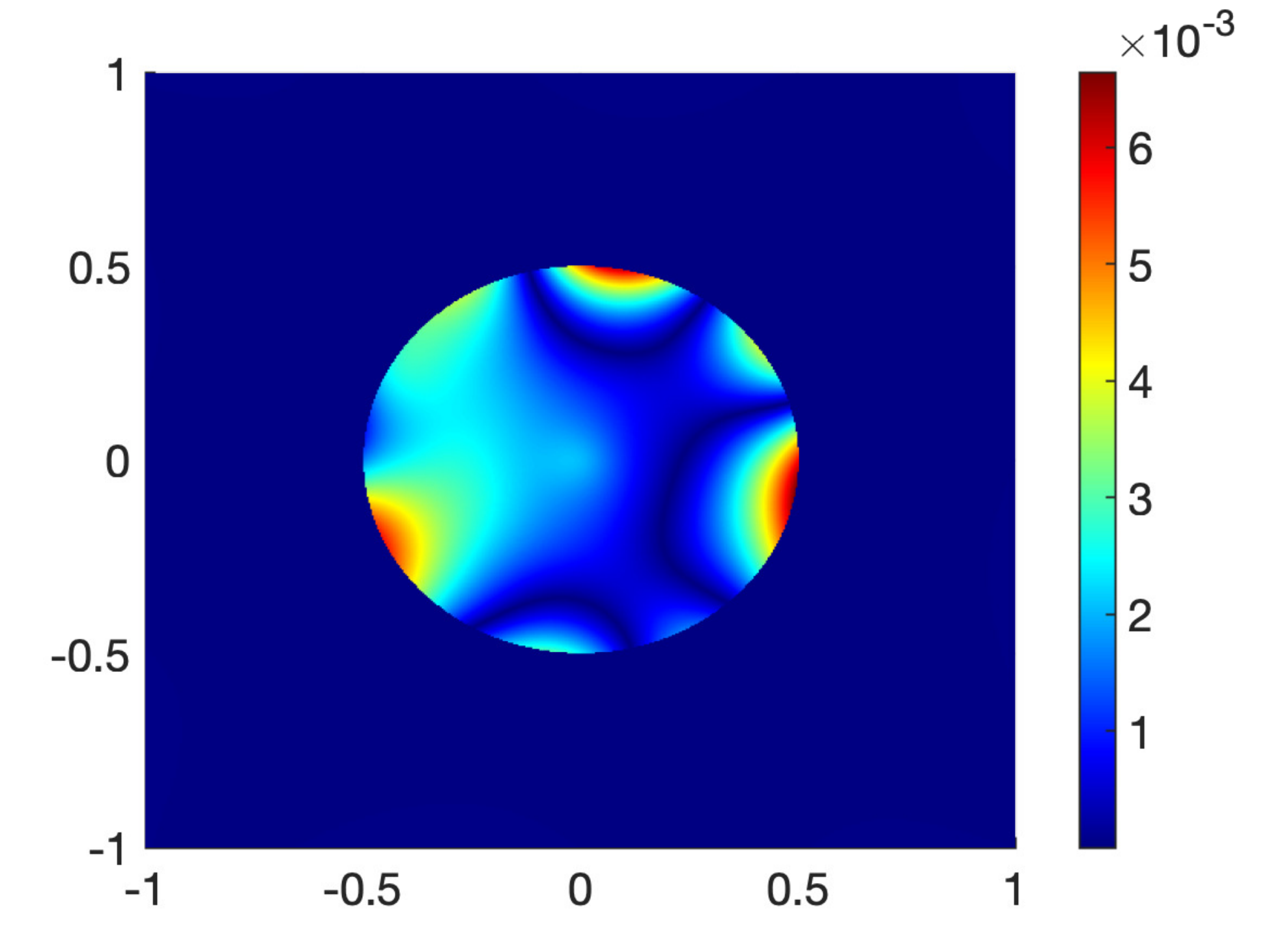}
		&
		\includegraphics[width=.32\textwidth]{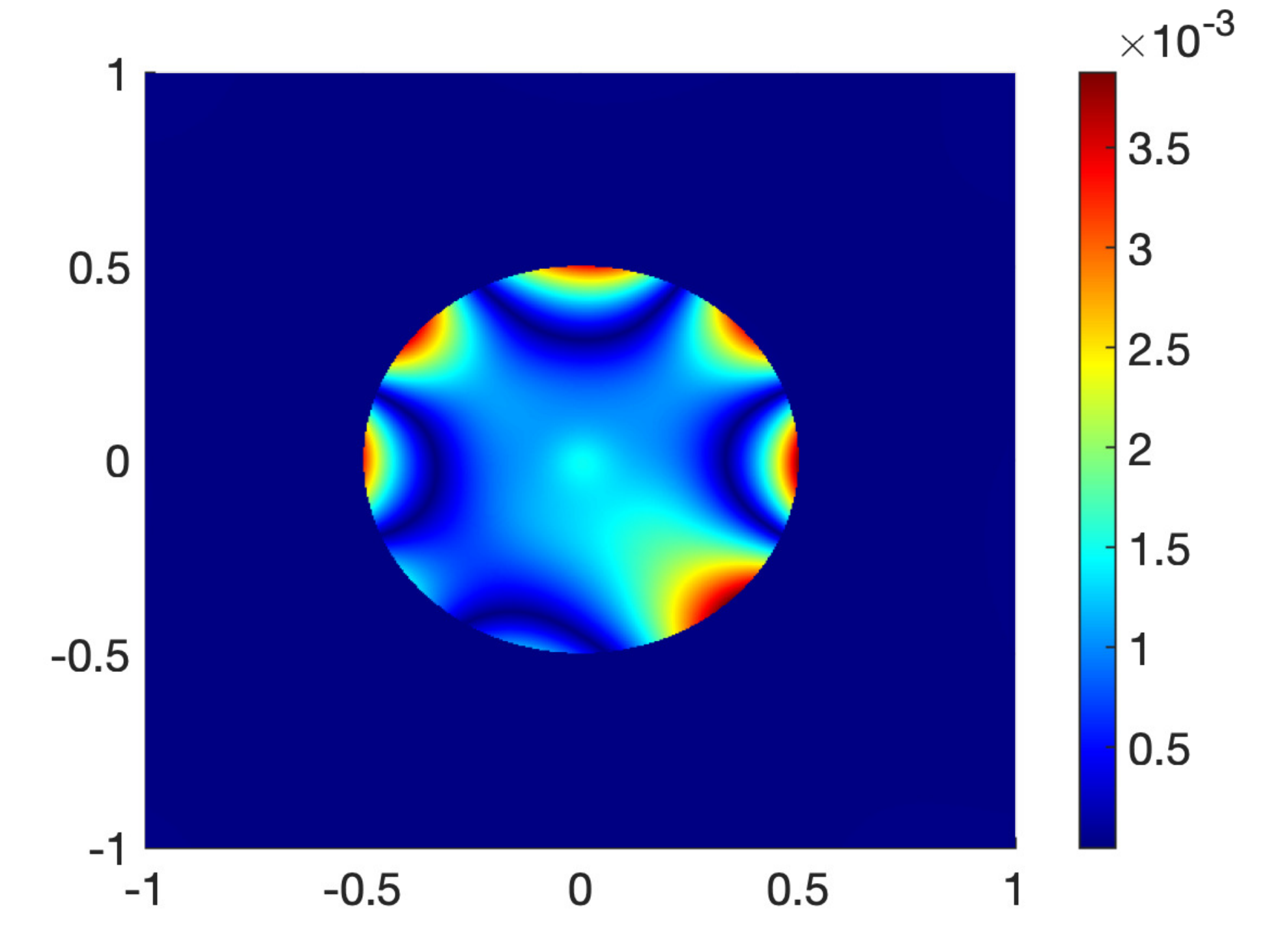}
		&
		\includegraphics[width=.32\textwidth]{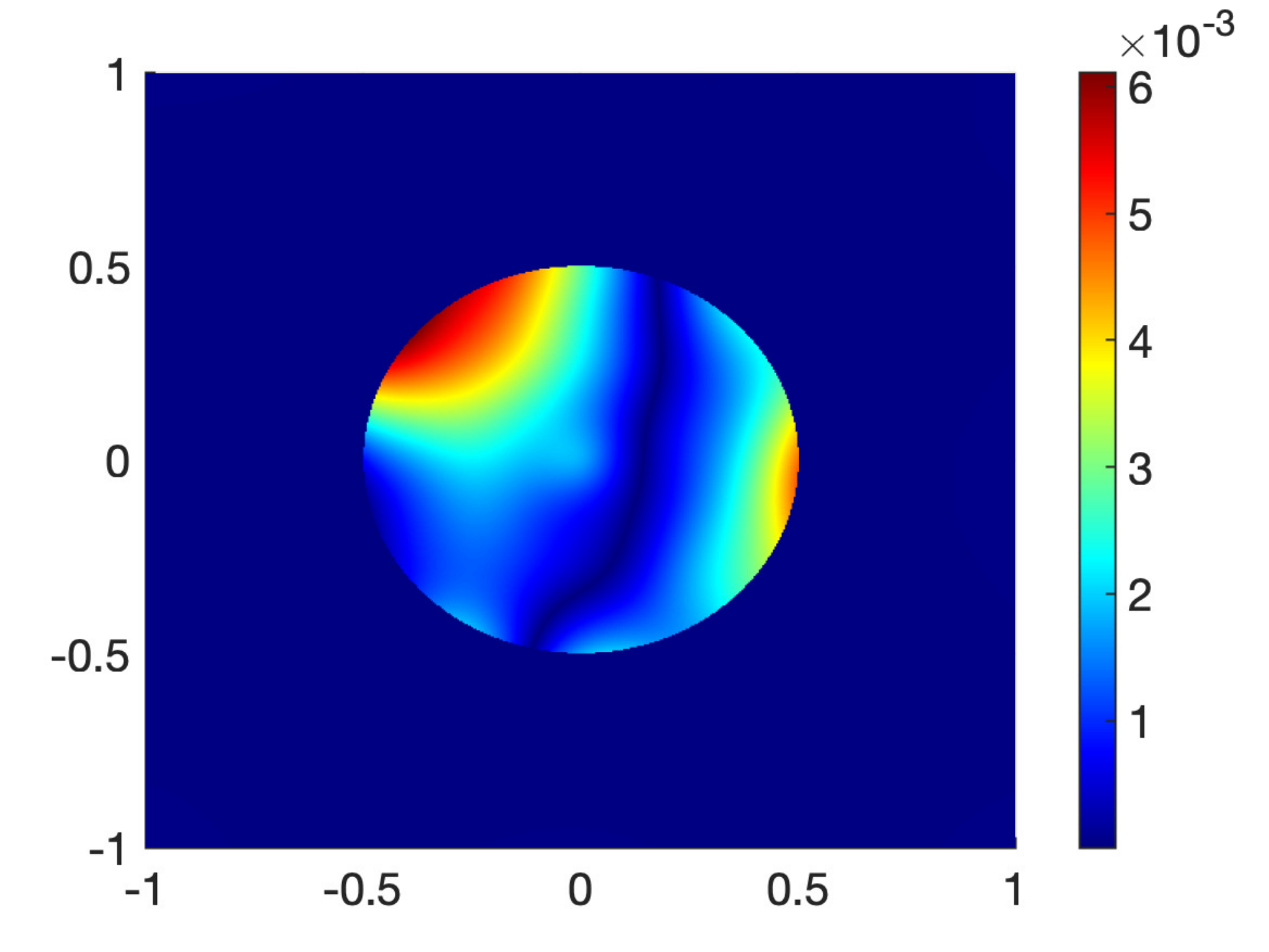}
		\\
		(a) & (b) & (c)
	\end{tabular}
	\caption{Example~\ref{Sect:NumTest1}. Error profiles of $\beta_1 = 1000$, $\beta_2 = 1$  on sampling size
		(a) $16 \times 16$;  (b) $32 \times 32$; (c) $64 \times 64$.
	}\label{Fig:Test1-a-tanh1}
\end{figure}

Figures \ref{Fig:Test1-1}(a) and \ref{Fig:Test1-1}(b) show the exact solutions when $\beta_1 =1000, \beta_2=1$ and $\beta_1 =1, \beta_2=1000$, respectively. We note that the solutions exhibit singularities (large jump in the derivative) across the interface. Moreover, Figure~\ref{Fig:Test1-a-tanh1} shows that the error profiles on the uniformly sampled points of sizes $16\times 16$, $32 \times 32$ and $64 \times 64$ with $10$ hidden layers for the case $\beta_1=1000$ and $\beta_2=1$. It is obvious that the errors are dominant near the interface. Based on this, it is then natural to sample the points adaptively based on the error. Inspired by the standard adaptive least-squares methods, we use the computable residual error, i.e, $\mathcal{L}_{\text{total}}$ ~\eqref{eqn:total-loss}, as the a posterior error indicator and investigate adaptive sampling techniques for effectively handling the solutions with singularities.  

\subsubsection{Case with $\beta_1 = 1000$ and $\beta_2 = 1$ }\label{Sect:NumTest1C}
In this example, we compare the approximations based on uniform sampling and adaptive sampling strategies.  The adaptive sampled points are obtained based on the residual error. More precisely, we start with uniformly sampled $10\times 10$ points and solve the interface problem.  Then we uniformly sample more points (in our experiments, we sample $5$ times more points in each direction for the next level) and compute the error indicator, i.e., the loss function~$\mathcal{L}_{\text{total}}$~\eqref{eqn:total-loss}, on those points.  Finally, we ranked the points according to the error indicator and add those ranked top $10\%$ to form the next adaptively refined level. This procedure is then repeated to generate more adaptively refined levels. 

\begin{figure}[H]
\centering
\begin{tabular}{ccc}
  \includegraphics[width=.31\textwidth]{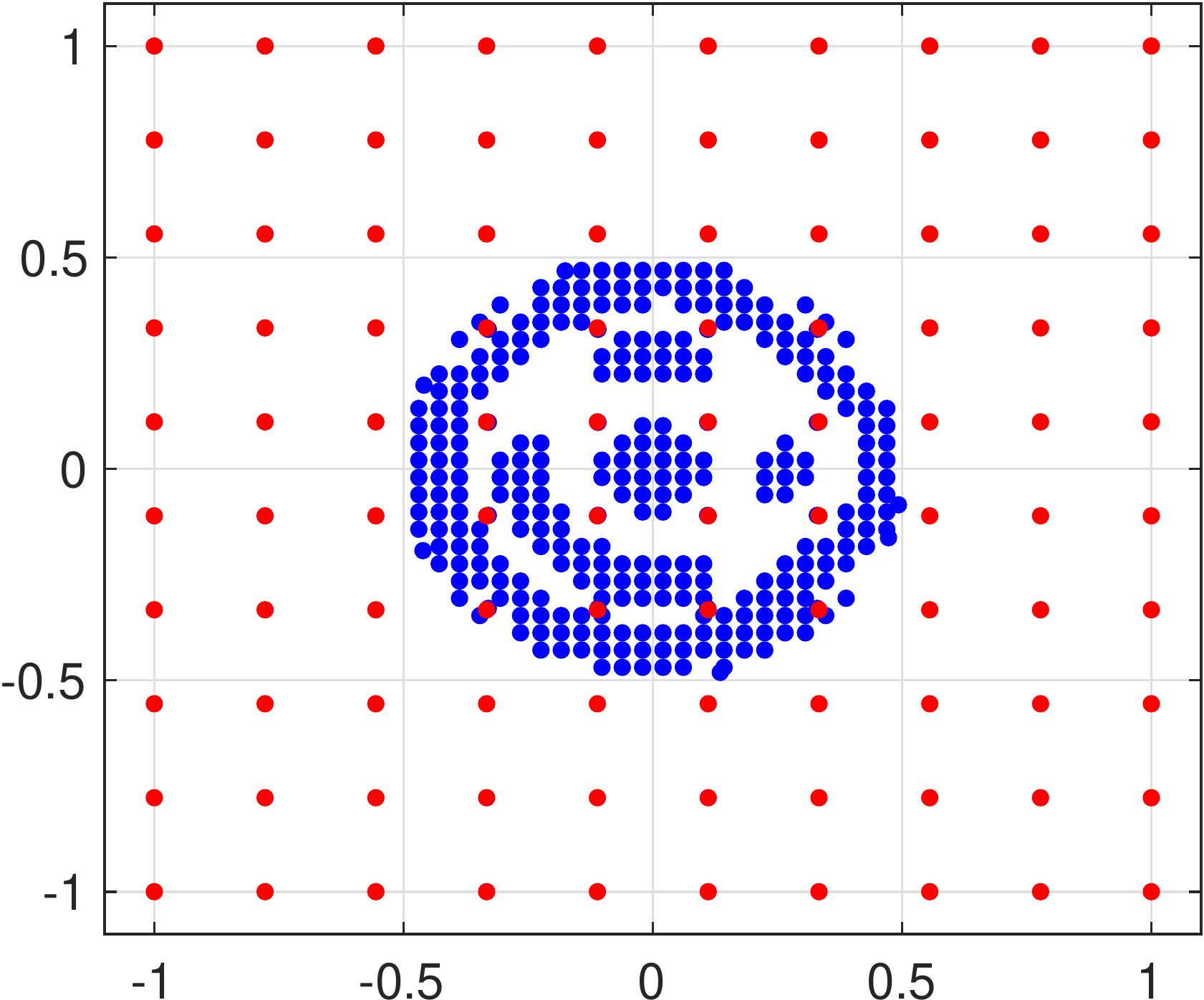}
&\includegraphics[width=.31\textwidth]{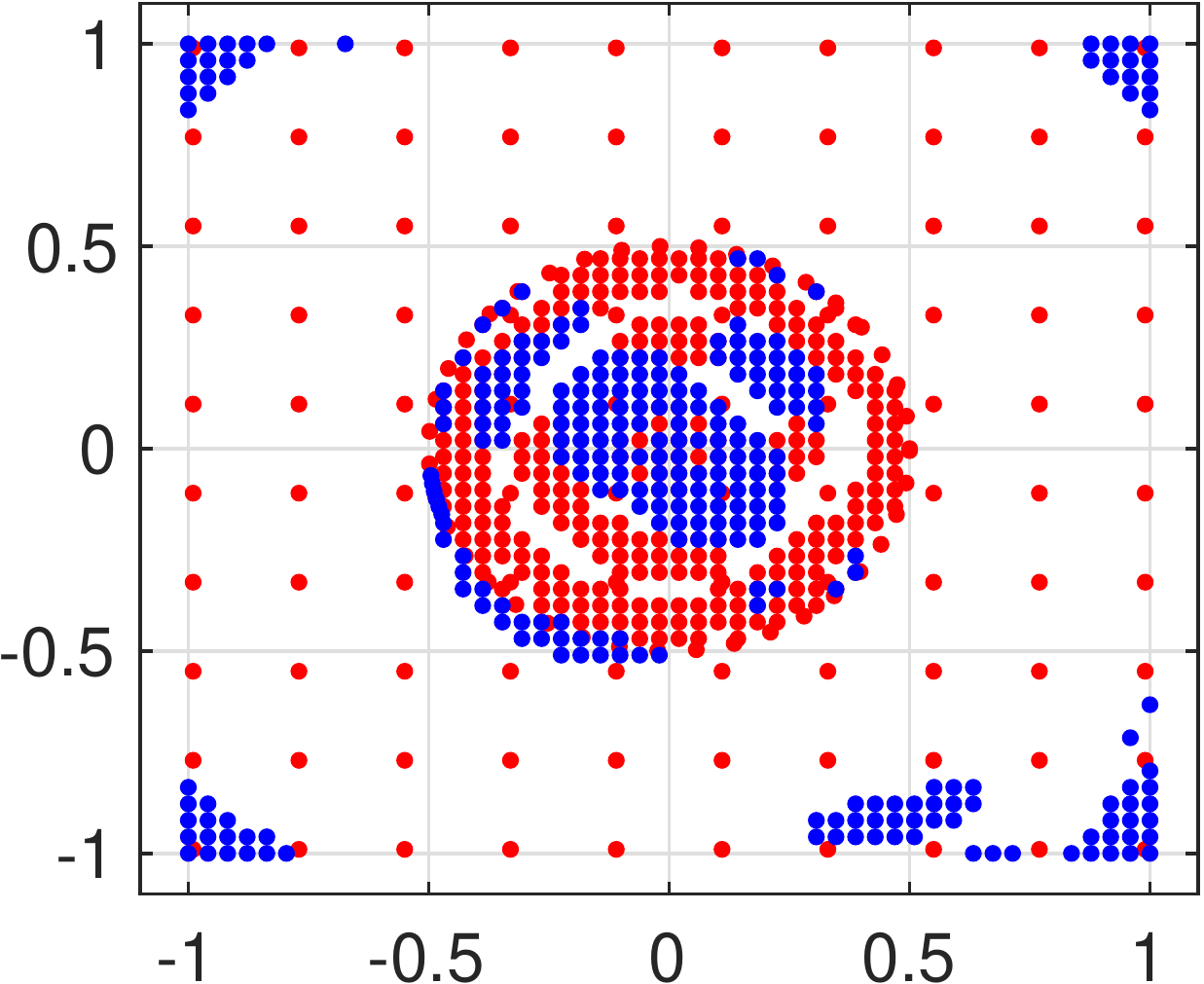}
&\includegraphics[width=.31\textwidth]{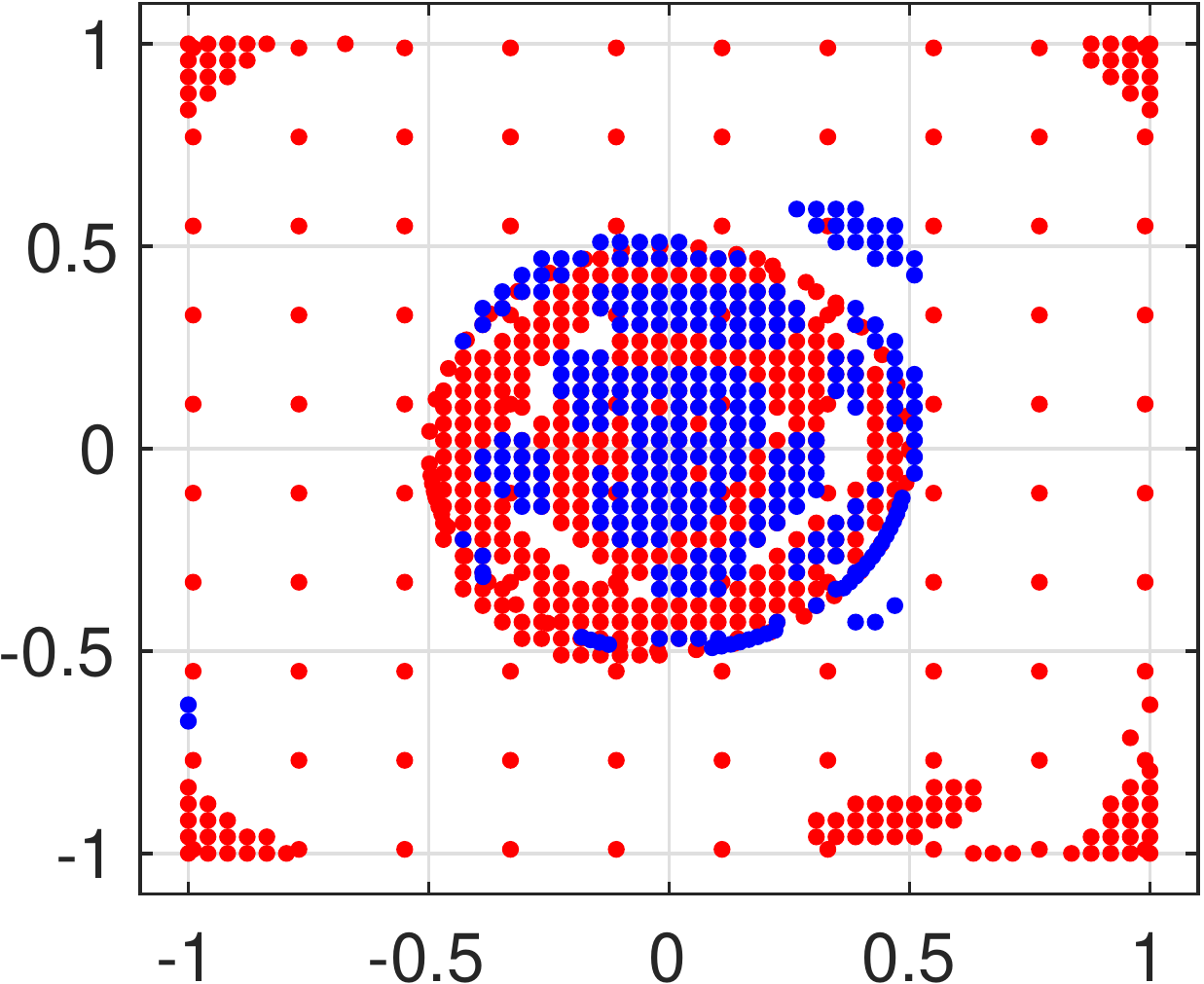}\\
(a) & (b) & (c)
\end{tabular}
\caption{ Example~\ref{Sect:NumTest1C} with $\beta_1 = 1000$ and $\beta_2 = 1$. Refinement (a) level 1; (b) level 2; (c) level 3. }\label{Fig:Test1C-1}
\end{figure}

Figure~\ref{Fig:Test1C-1} plots the sampling points for the first three refinements as the blue dots denoting the newly added points and red dots denoting the existing sampling points. It can be seen that relatively more points are added near the interface, which captures the singularities of the solution. 

In Table~\ref{Tab:Test1C-1}, we quantify and compare the performances
in terms of the relative error in the $L^2$ norm between the uniform sampling and adaptive sampling strategies. Two different number of hidden layers, i.e., $4$ and $6$, are used in the DNN structure.
As we can see, the errors decrease effectively as we adaptive sample points. Furthermore, with only three refinements, the errors on the third adaptive refinements are comparable to that of $50 \times 50$ uniformly sampled points. Note that the number of points in the adaptive setting is about $60\%-70\%$ less than that of the uniform case.

\begin{table}[H]
\caption{Example~\ref{Sect:NumTest1C}.
{$\dfrac{\|u - \mathcal{U_{NN}}\|_\Omega}{\|u\|_\Omega}$ with $\beta_1 = 1000$ and $\beta_2 = 1$.}
}
\label{Tab:Test1C-1}
\centering
\begin{tabular}{c||cccc||cccc}
\hline\hline
&\multicolumn{4}{c||}{4 Layers} &\multicolumn{4}{c}{6 Layers}\\
Grids & $M_{1}$ & $M_{2}$  & $M_{\Gamma}$ &Error & $M_{1}$ & $M_{2}$  & $M_{\Gamma}$  &6 Layers \\[5pt] \hline
Uniform $10\times 10$ &16 &84 &32 & 5.4484e-03 &16 &84 &32 &1.3514e-02\\
Refined Level 1 &304 &86   &32 &1.4922e-03 &303 &87&32 &4.4854e-03\\
Refined Level 2 &405 &202 &37 &1.1534e-03 &375 &266 &47&9.4701e-04\\
Refined Level 3 &450 &251 &59 &1.0530e-03 &406 &477 &53&8.5256e-04\\ \hline
Uniform $50\times 50$ &484     &2016  &160 & 2.4664e-03 &484     &2016  &160 &8.3150e-04\\
\hline\hline
\end{tabular}
\end{table}

\subsubsection{Case with $\beta_1 = 1$ and $\beta_2 = 1000$ }\label{Sect:NumTest1D}
\begin{figure}[H]
\centering
\begin{tabular}{ccc}
  \includegraphics[width=.31\textwidth]{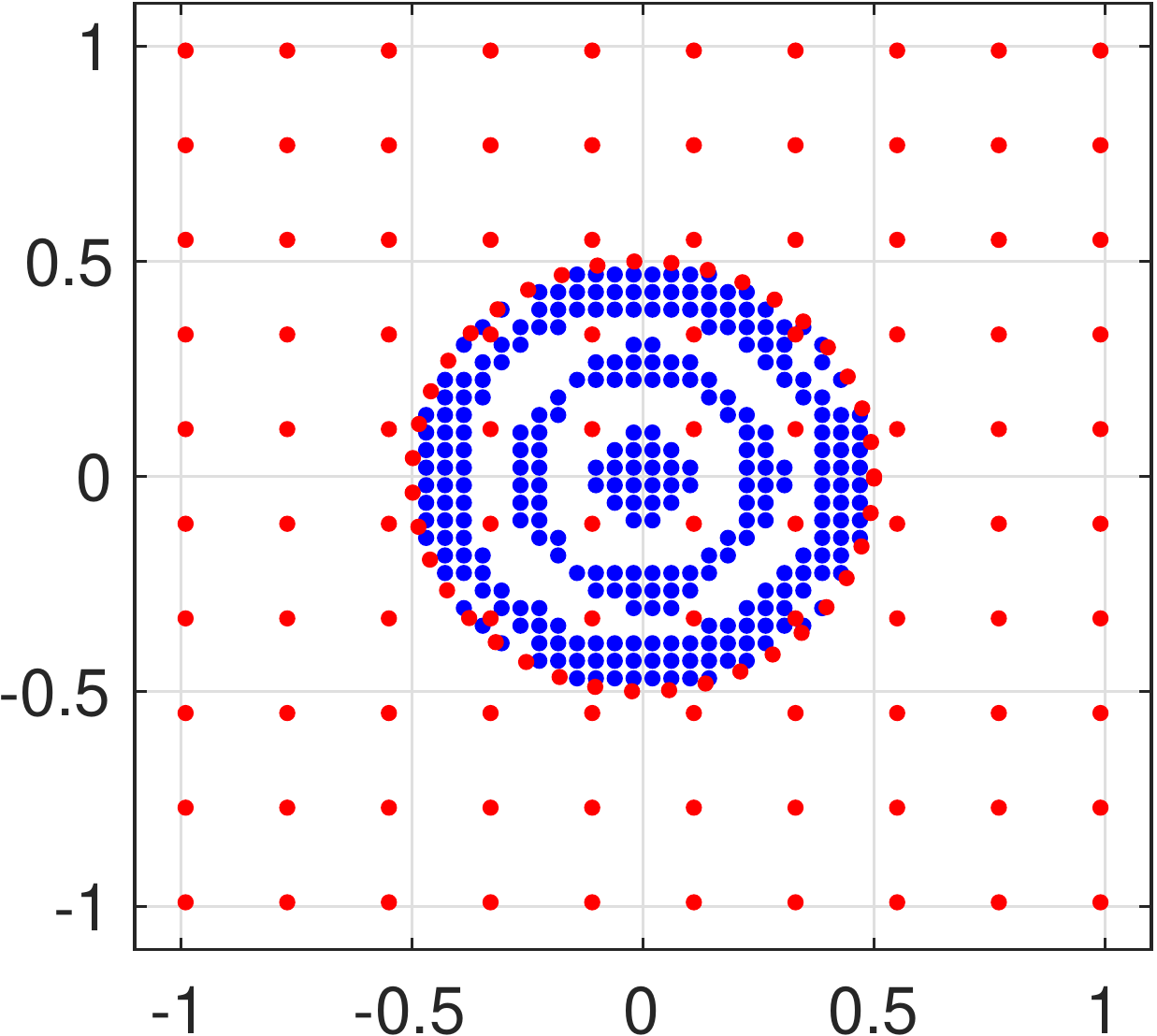}
&\includegraphics[width=.31\textwidth]{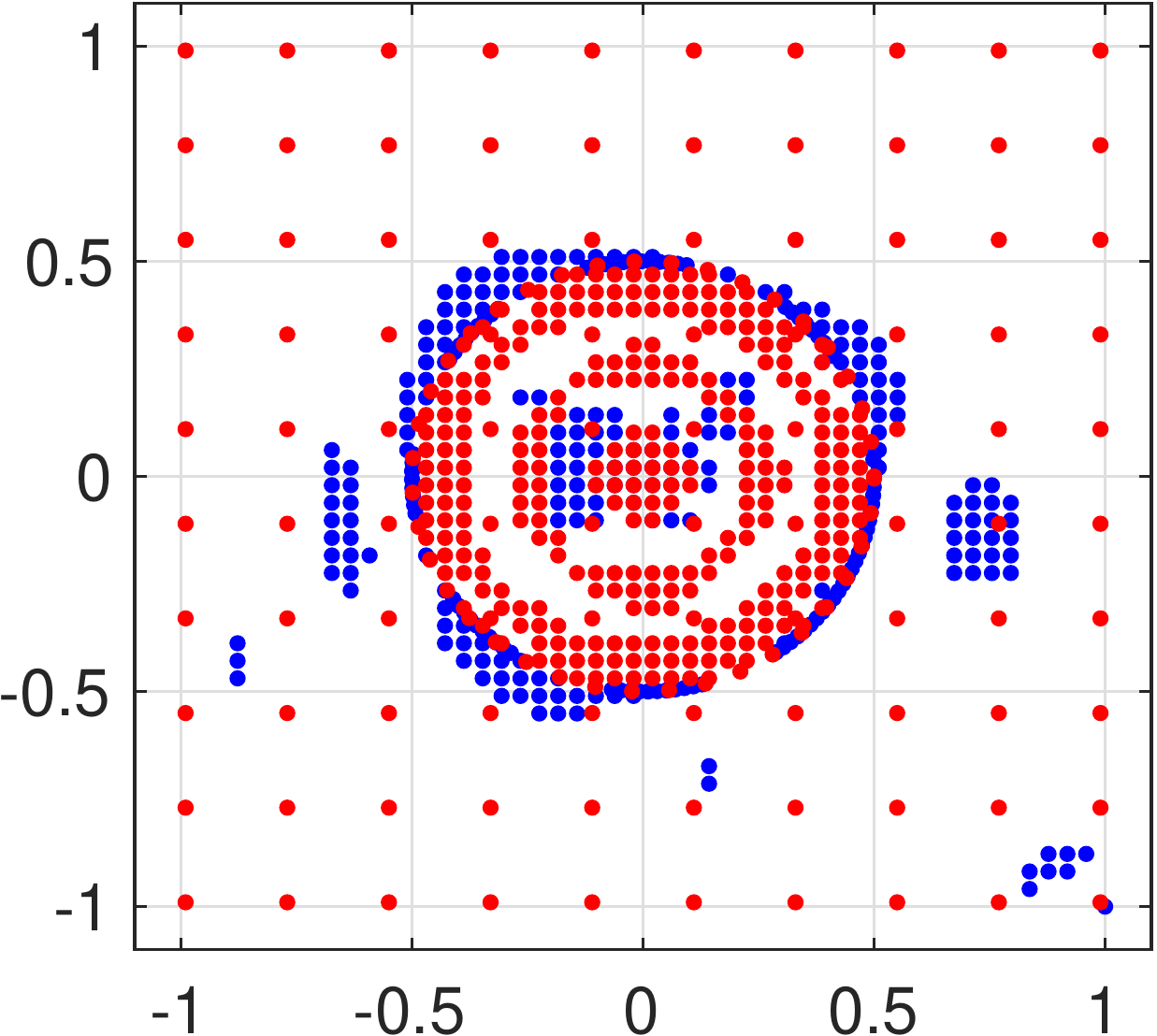}
&\includegraphics[width=.31\textwidth]{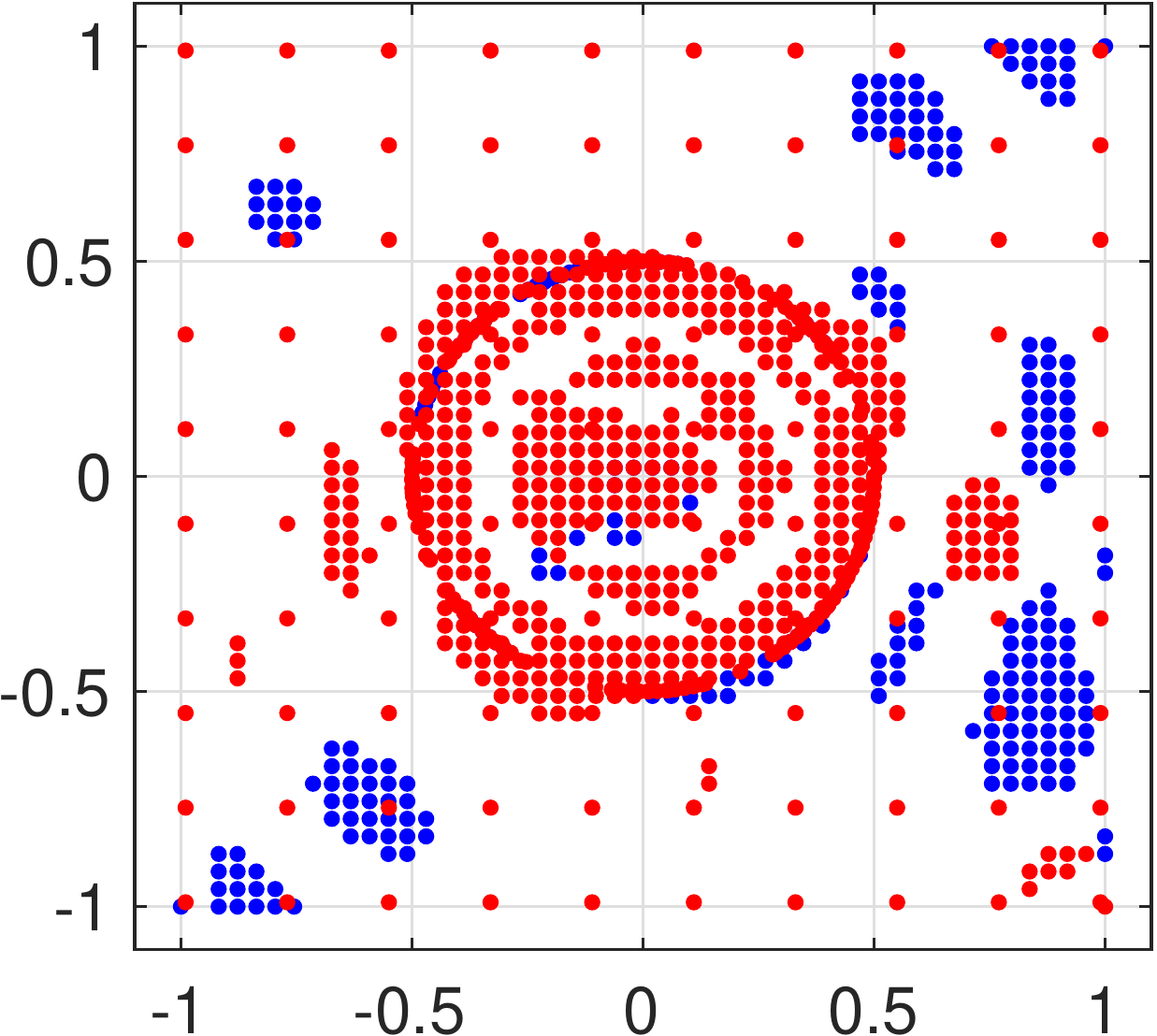}\\
(a) & (b) & (c)
\end{tabular}
\caption{ Example~\ref{Sect:NumTest1D} with $\beta_1 = 1$ and $\beta_2 = 1000$. Refinement (a) level 1;(b) level 2; (c) level 3. }
\label{Fig:Test1D-1}
\end{figure}

We now consider the case that $\beta_1=1$ and $\beta_2=1000$.  The adaptive sampled points for the first three refinements are given in Figure \ref{Fig:Test1D-1} and the comparison of the errors is listed in Table \ref{Tab:Test1D-1}. Similar conclusions as in Section~\ref{Sect:NumTest1C} can be drawn.

\begin{table}[H]
\caption{Example~\ref{Sect:NumTest1D}.
{$\dfrac{\|u - \mathcal{U_{NN}}\|_\Omega}{\|u\|_\Omega}$ with $\beta_1 = 1$ and $\beta_2 = 1000$.}
}\label{Tab:Test1D-1}
\centering
\begin{tabular}{c||cccc||cccc}
\hline\hline
&\multicolumn{4}{c||}{4 Layers}&\multicolumn{4}{c}{6 Layers}\\
Grids & $M_{1}$ & $M_{2}$  & $M_{\Gamma}$ &Error & $M_{1}$ & $M_{2}$  & $M_{\Gamma}$  &Error  \\[5pt] \hline
Uniform $10\times 10$ &16 &84       &32    &5.0393e-02 &16 &84       &32 &1.4533e-02\\
Refined Level 1 &304 &86   &32    &2.7827e-02 &304 &86 &32&6.1598e-03\\
Refined Level 2 &414 &119 &106  &3.2294e-03 &352 &232 &102&3.8287e-03\\
Refined Level 3 &440 &266 &130 &3.1326e-03 &361 &449 &115&1.6241e-03\\ \hline
Uniform $50\times 50$ &484 &2016 &160 &3.1357e-03 &484 &2016 &160 &1.6239e-03\\
\hline\hline
\end{tabular}
\end{table}

\subsection{Example 6. Flower Shape Interface}\label{Sect:NumTest2}
\begin{figure}[H]
\centering
\begin{tabular}{cc}
\includegraphics[width=.45\textwidth]{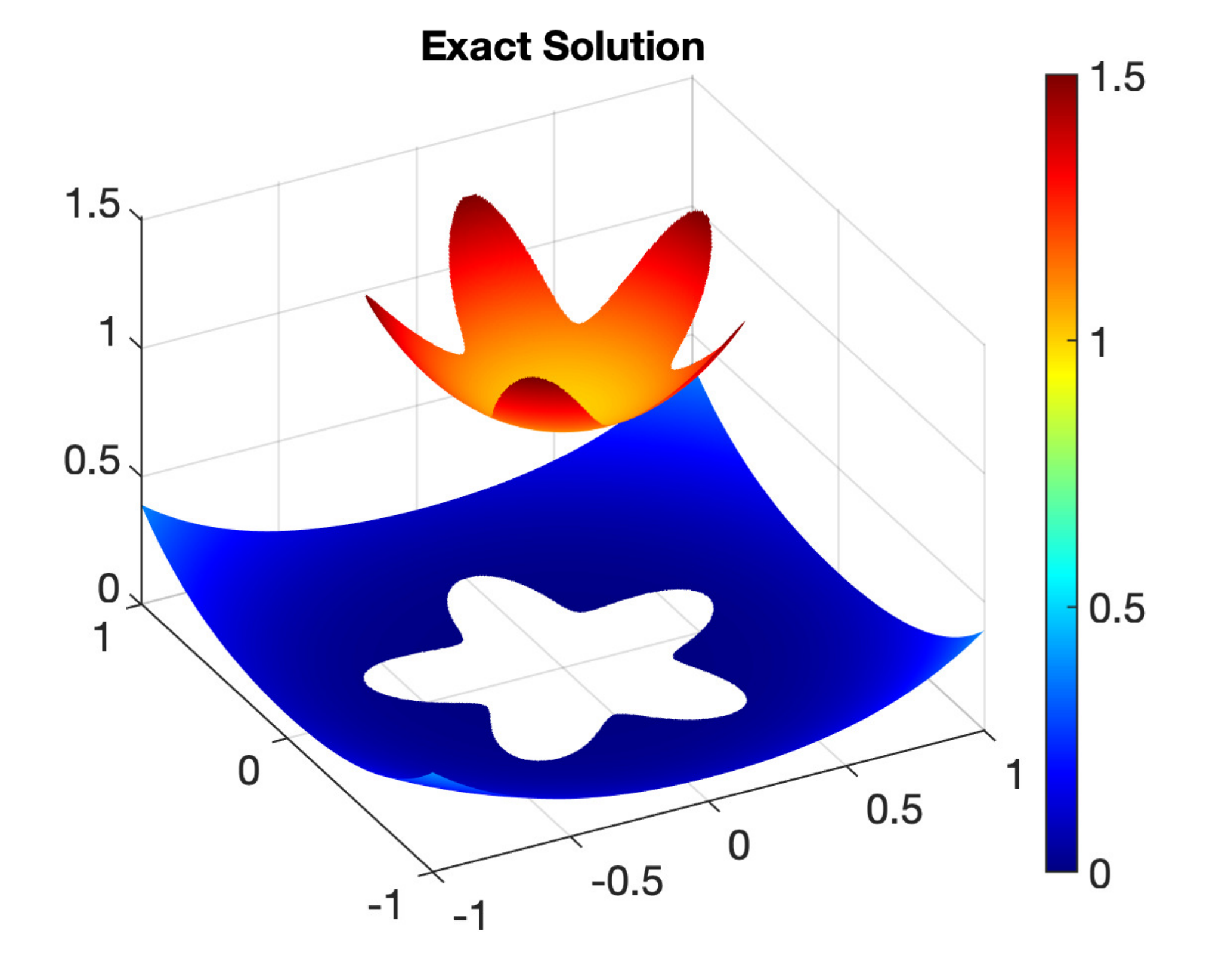}
&\includegraphics[width=.45\textwidth]{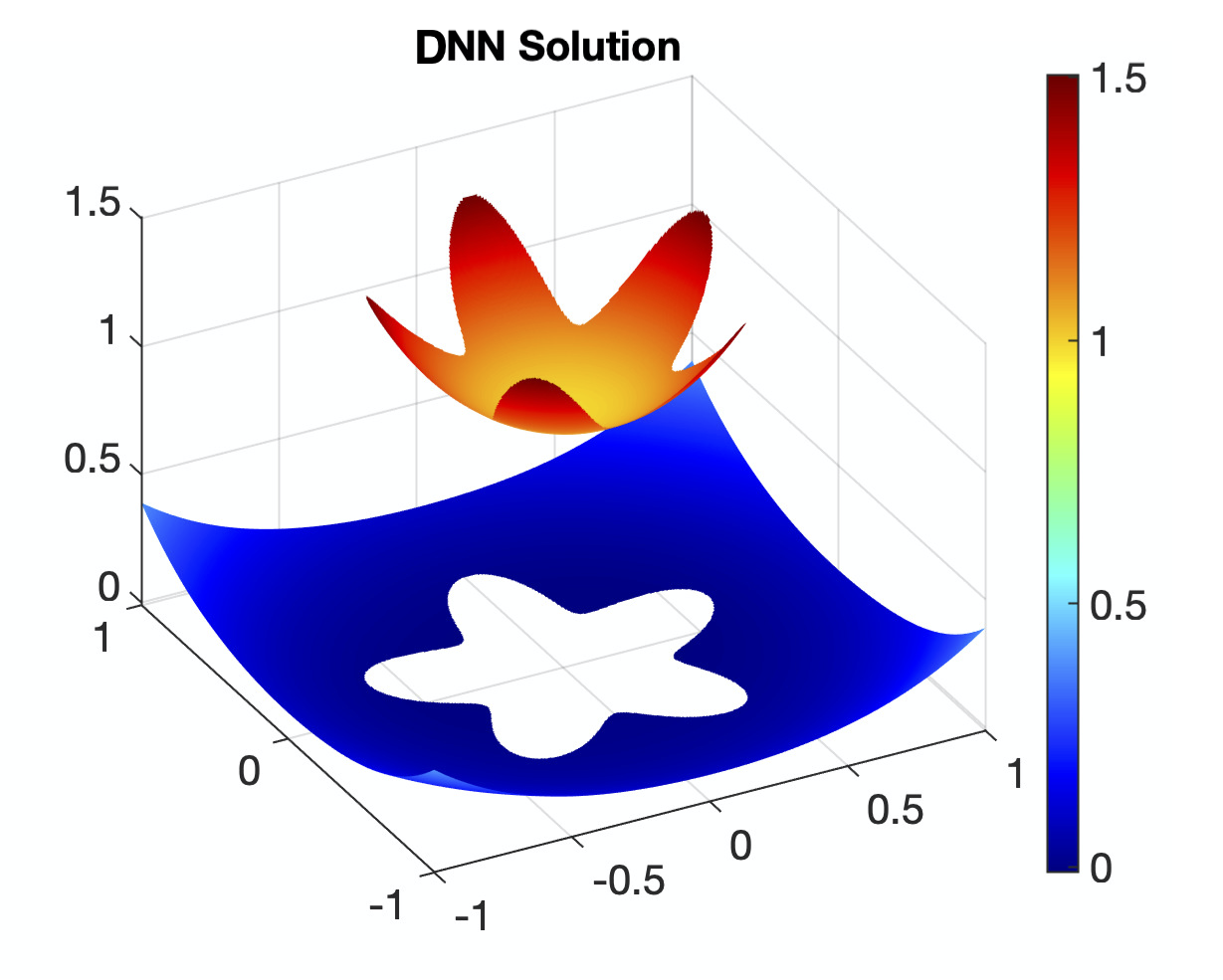}\\
(a) & (b) 
\end{tabular}
\caption{ Example~\ref{Sect:NumTest2}. (a) plot of exact solution; (c) plot of numerical solution. }\label{Fig:Test2-1}
\end{figure}

In this example, we consider a more complicated interface, i.e., the flower-shaped interface problem (\ref{eq:pde})-(\ref{eq:pde-bc}) with non-homogeneous jump condition to test our algorithm with adaptive sampling. The interface $\Gamma$ is given by the following equation,
$$
r=\frac{1}{2}+\frac{\sin(5\theta)}{7}.
$$ 
The exact solution is chosen as (see Figure~\ref{Fig:Test2-1}),
\begin{eqnarray}
u(\bm{x})=\begin{cases}
\exp(x^2+y^2),&\mbox{ if } \bm{x}\in\Omega_1,\\
0.1(x^2+y^2)^2-0.01\ln(2\sqrt{x^2+y^2}),&\mbox{ if } \bm{x} \in\Omega_2,
\end{cases}
\end{eqnarray}
and $\beta_1 = 10$ and $\beta_2 = 1$.
The jump conditions $g_j$ and $g_f$ are then computed by the exact solution and $\beta$. 
 Note that the coefficient contrast is mild in this case.

\begin{figure}[H]
\centering
\begin{tabular}{cccc}
\includegraphics[width=.21\textwidth]{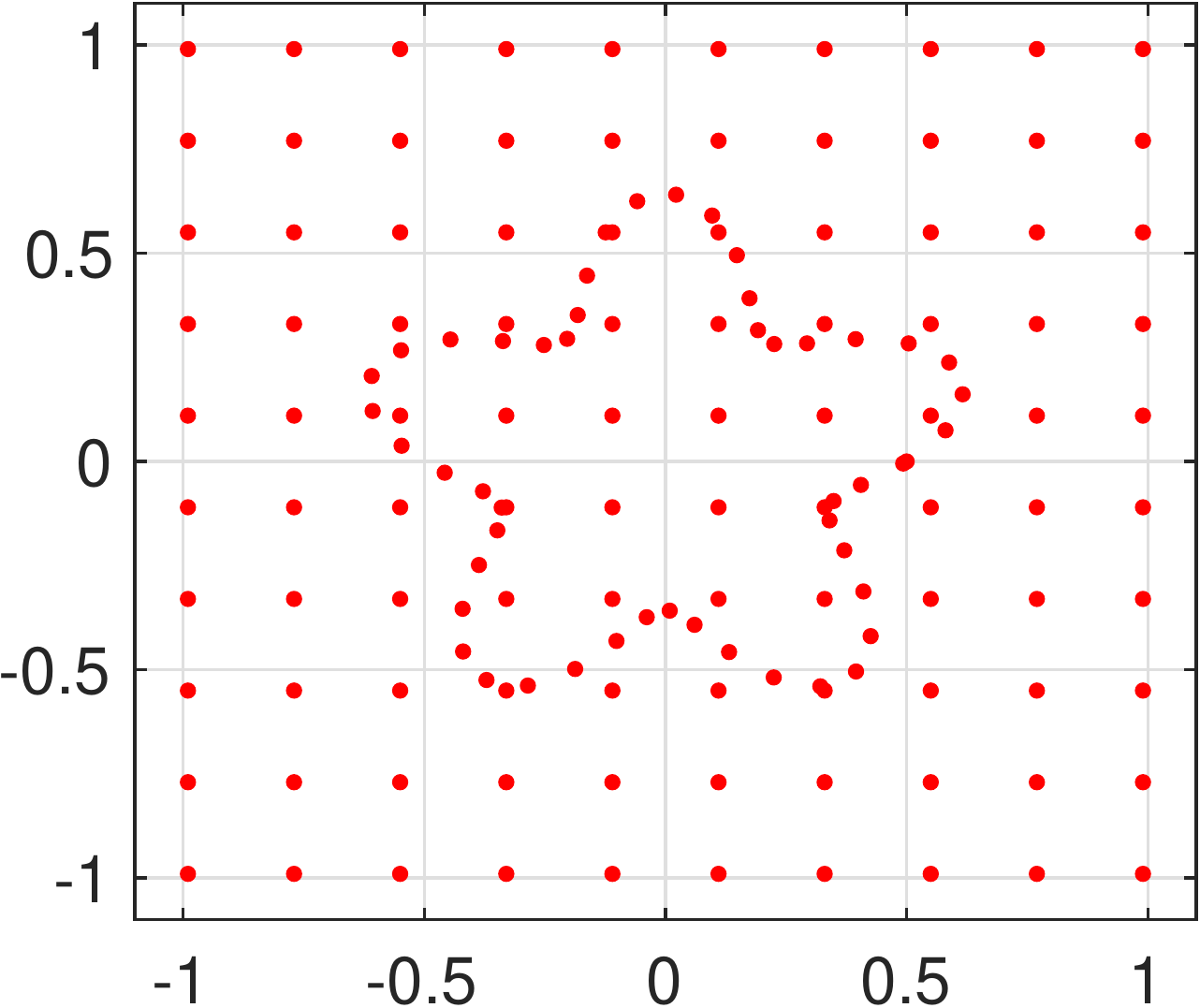}
&\includegraphics[width=.21\textwidth]{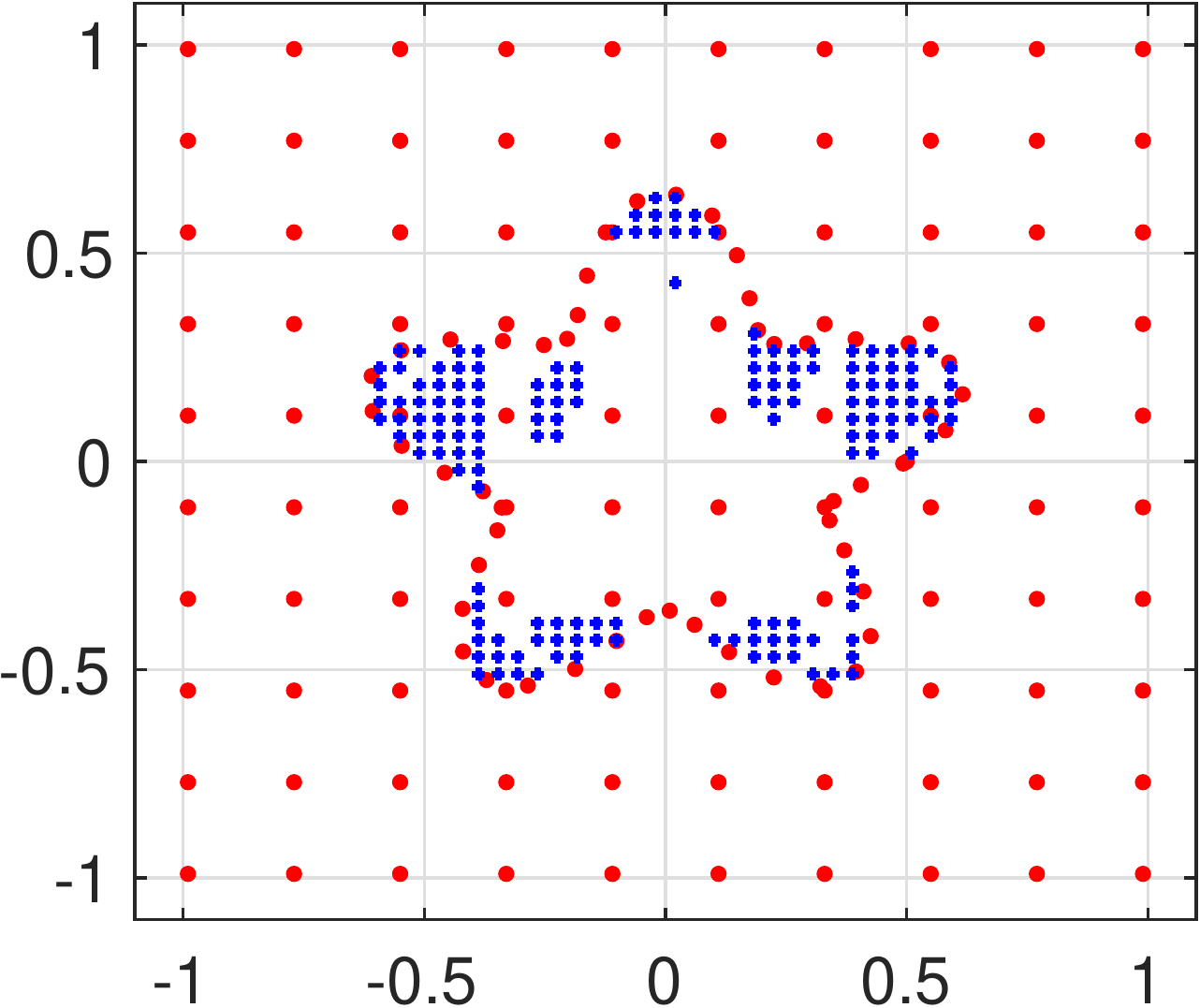}
&\includegraphics[width=.21\textwidth]{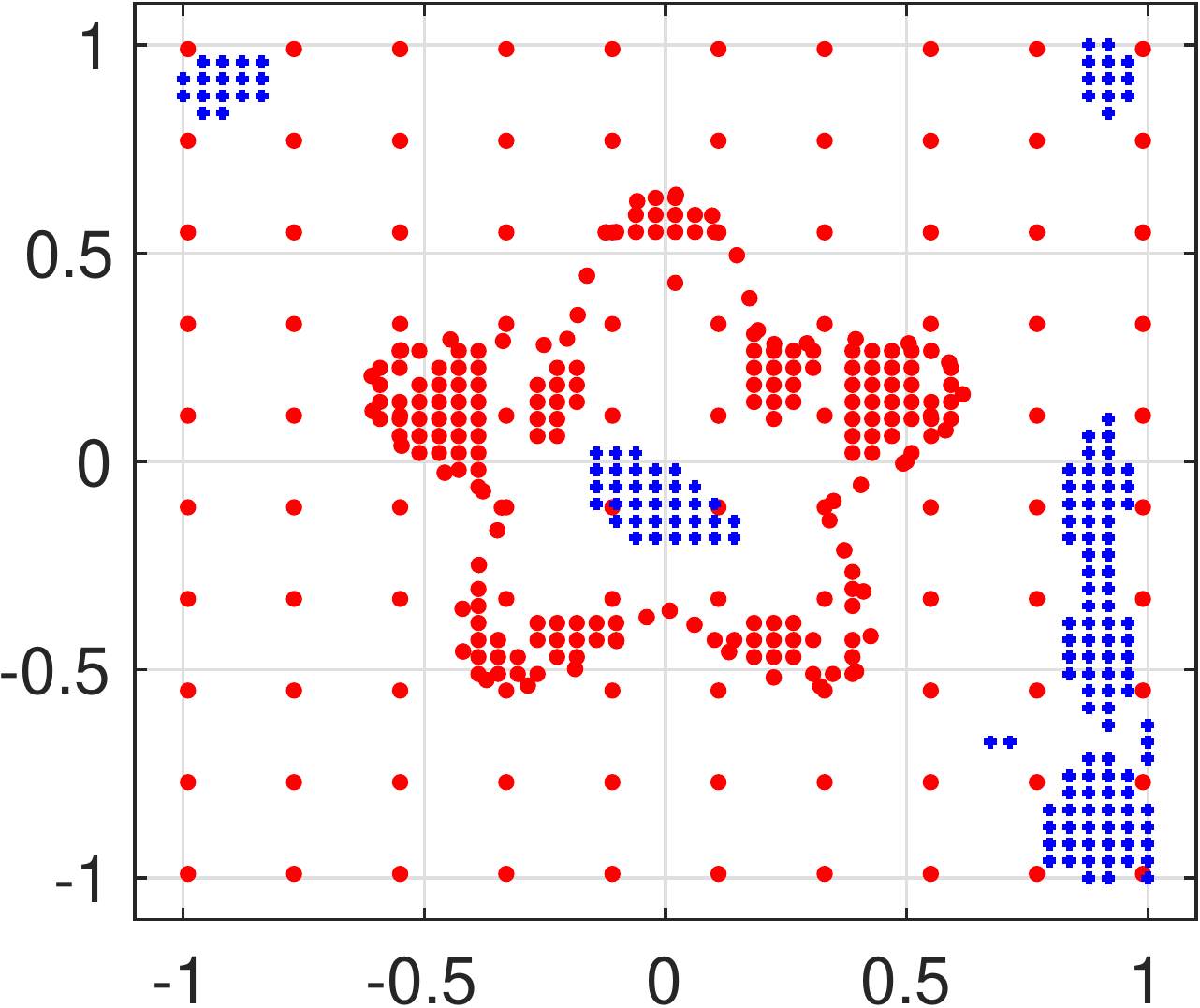}
&\includegraphics[width=.21\textwidth]{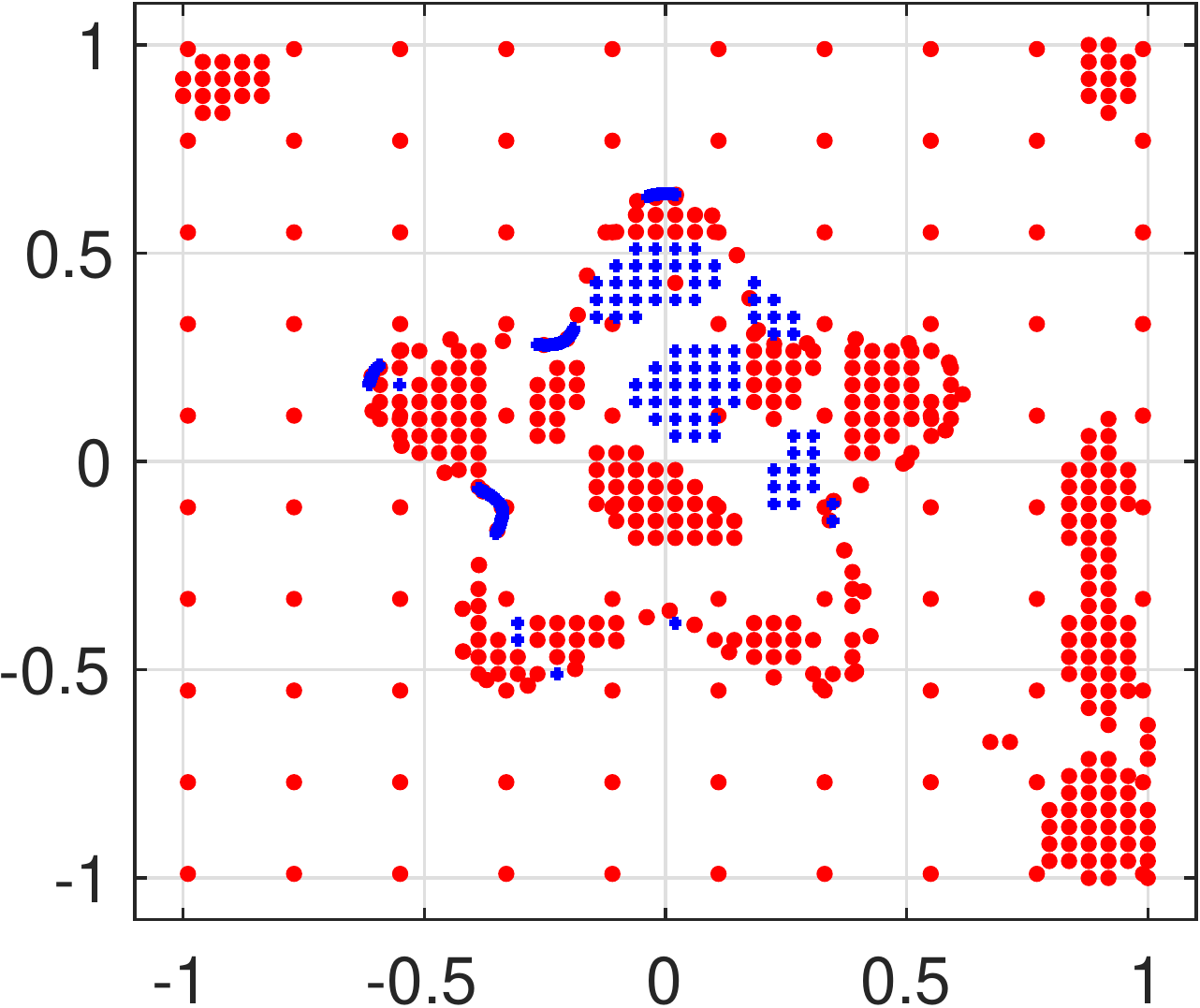}\\
(a) & (b) & (c) & (d)
\end{tabular}
\caption{ Example~\ref{Sect:NumTest2}. Refinement (a) level 0; (b) level 1; (c) level 2; (d) level 3. }
\label{Fig:Test2-3}
\end{figure}

\begin{figure}[H]
\centering
\begin{tabular}{cccc}
\includegraphics[width=.21\textwidth]{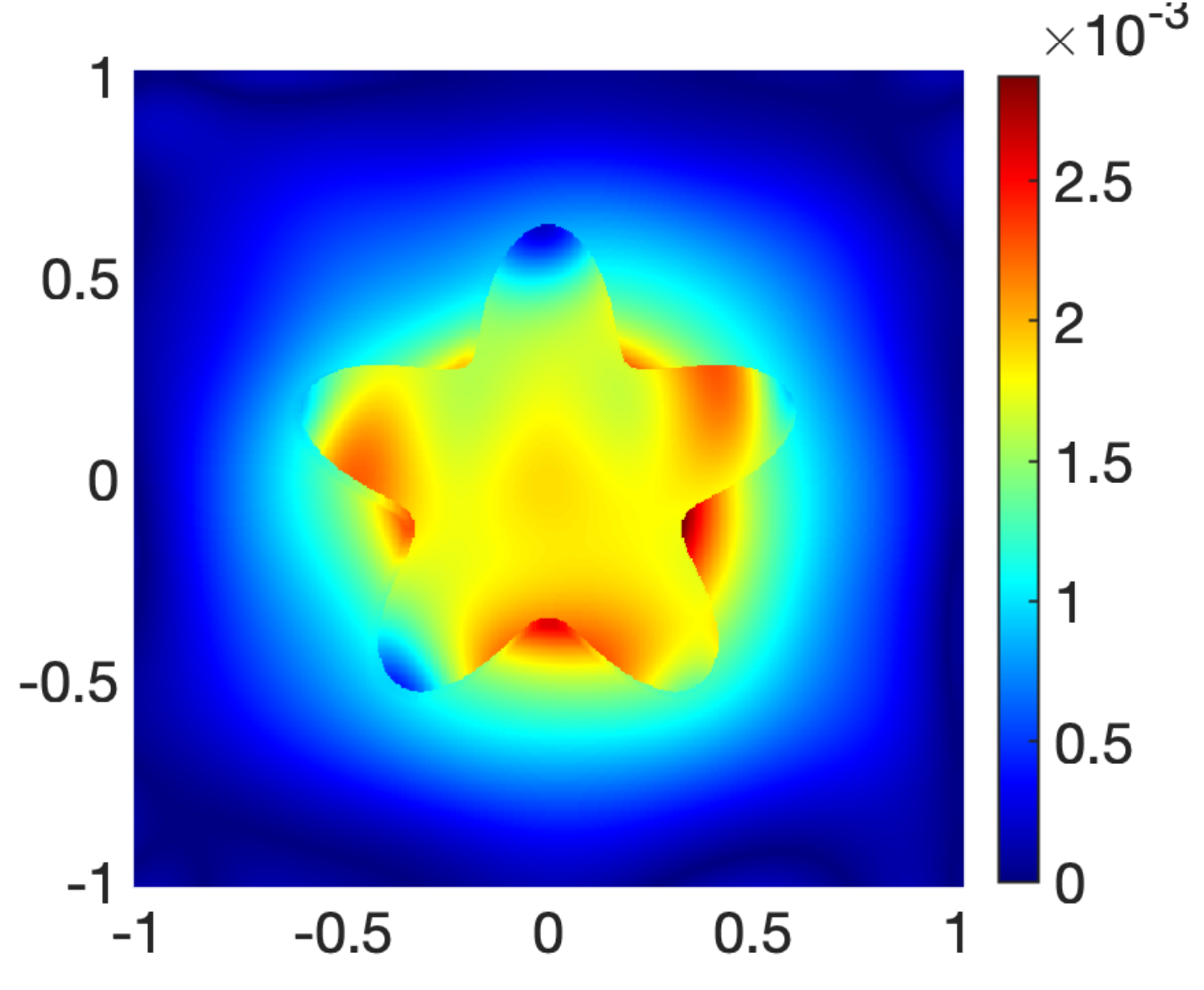}
&\includegraphics[width=.21\textwidth]{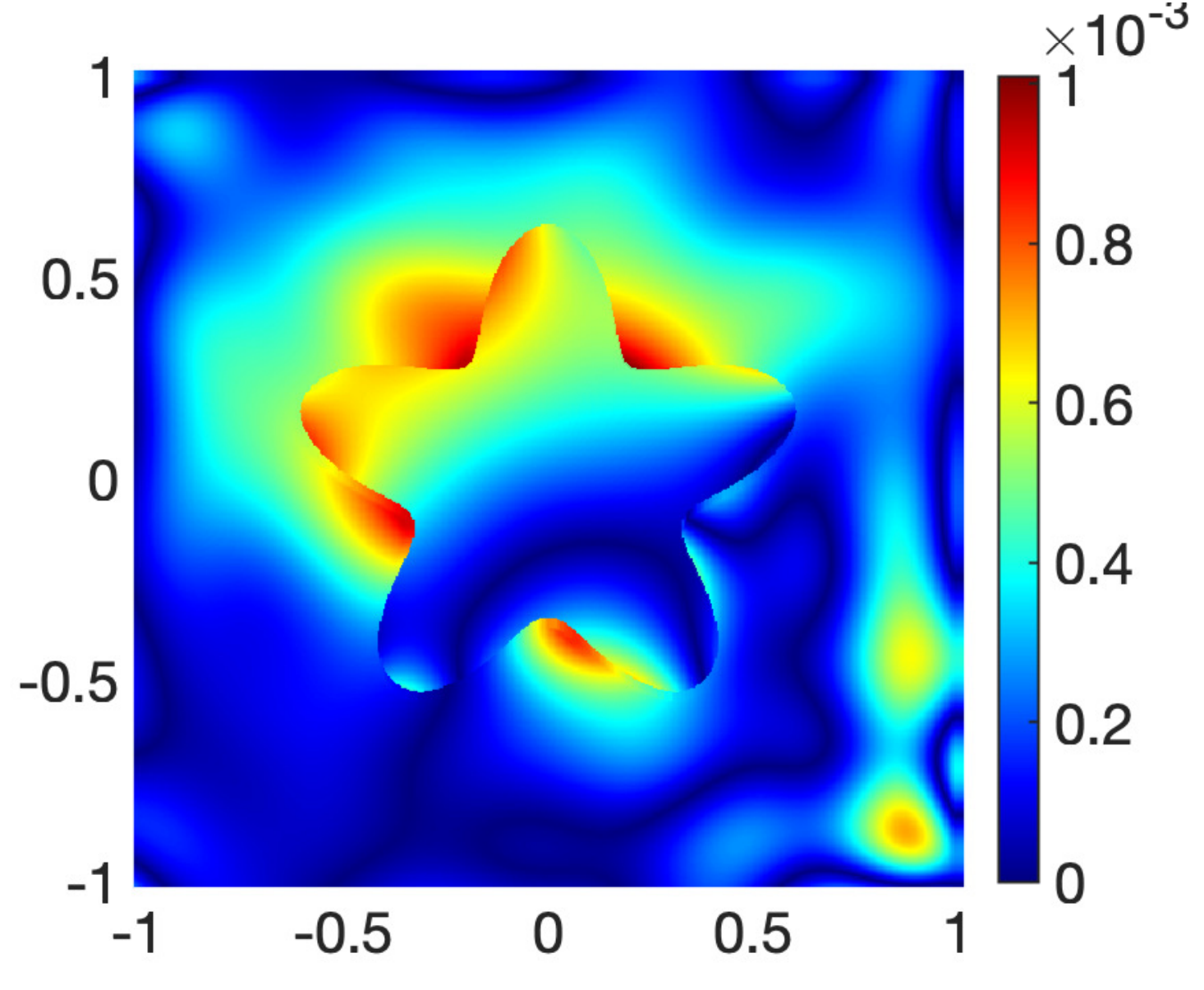}
&\includegraphics[width=.21\textwidth]{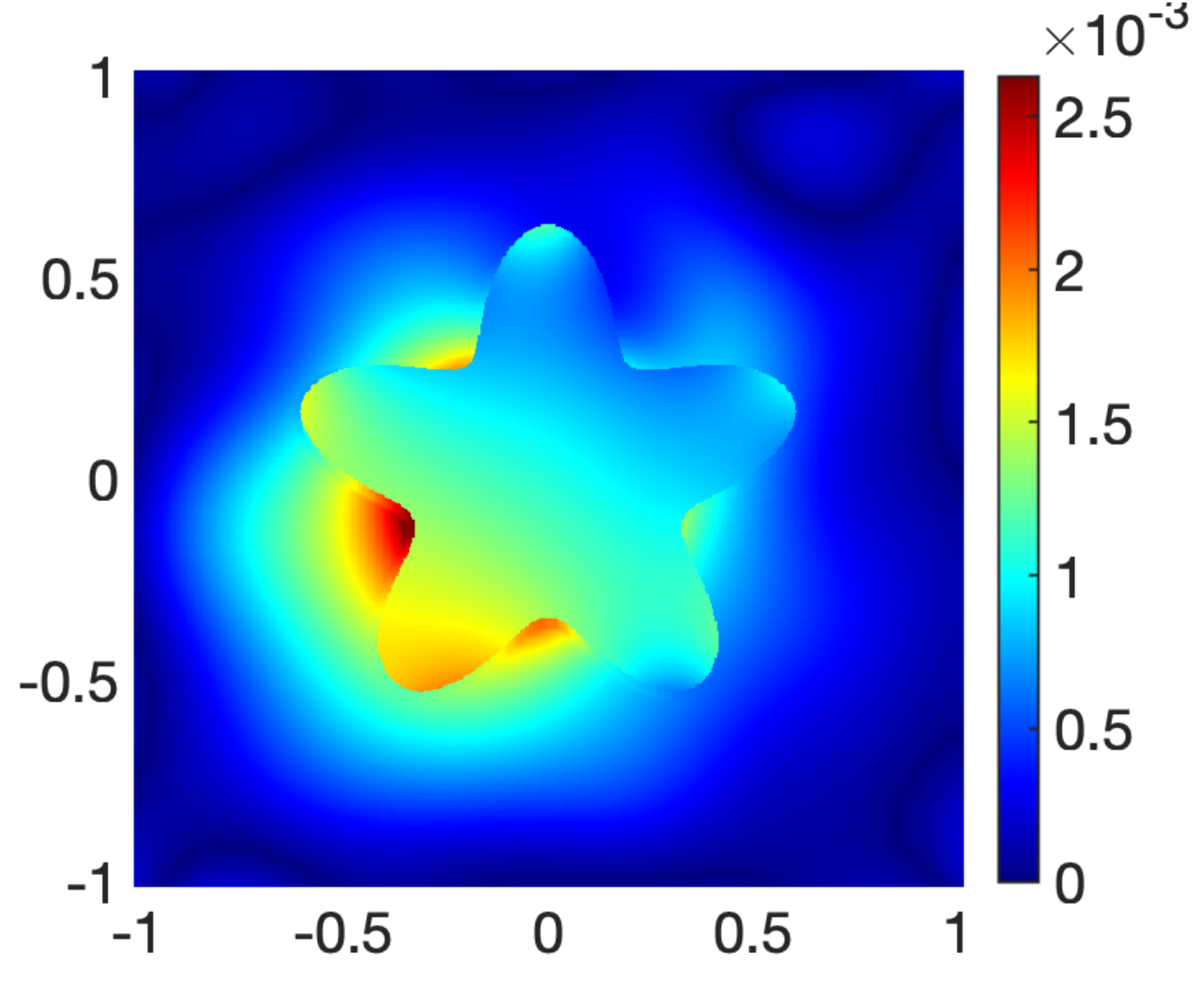}
&\includegraphics[width=.21\textwidth]{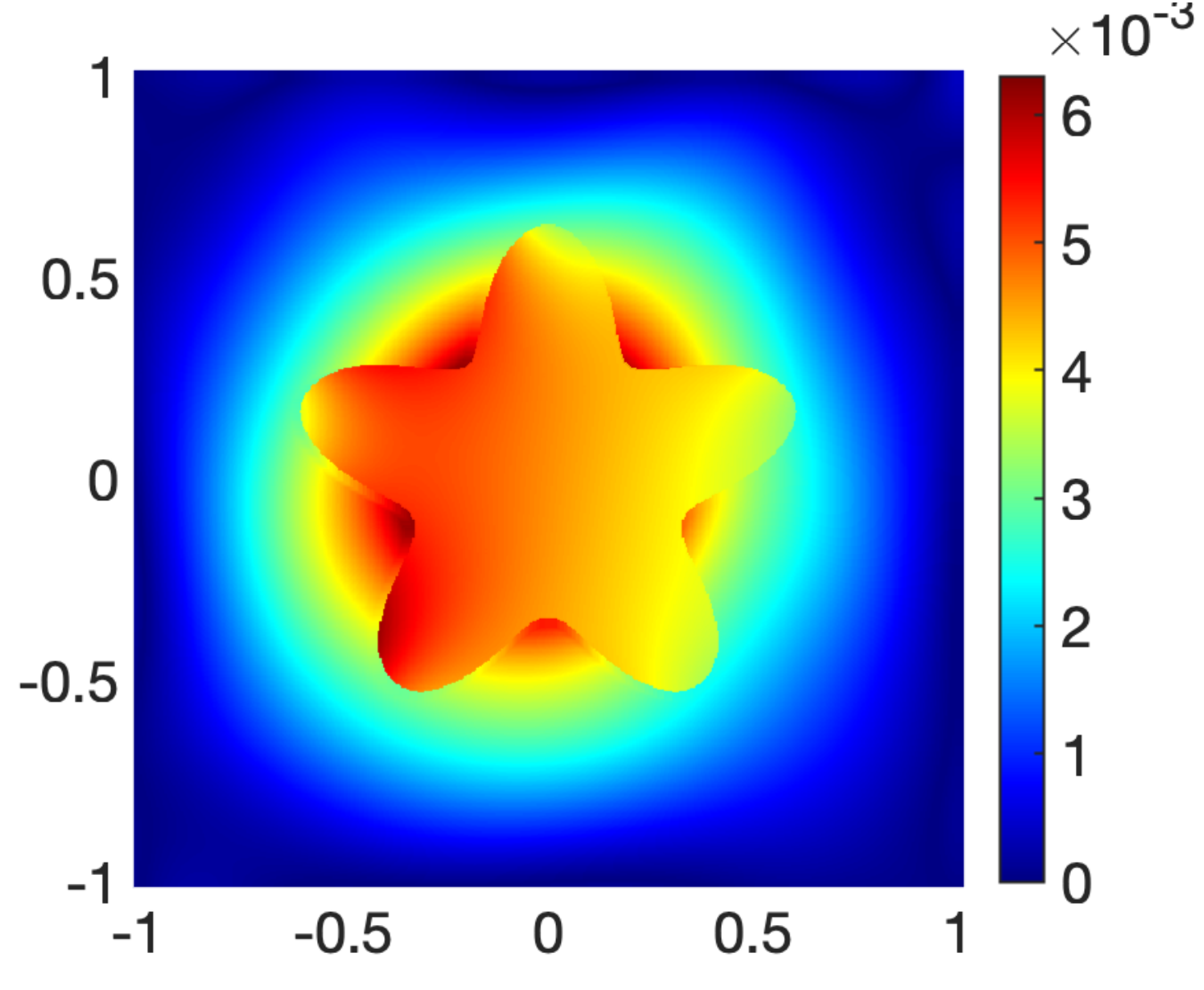}\\
(a) & (b) & (c) &(d)
\end{tabular}
\caption{ Example~\ref{Sect:NumTest2}. Plot of errors on (a) uniform $10\times 10$ grid; (b) adaptive refinement 1; (c) adaptive refinement 2; (d) adaptive refinement 3. }
\label{Fig:Test2-2}
\end{figure}

Figure \ref{Fig:Test2-3} provides the adaptive sampled points of the first three refinements. Again, red dots denoting the exiting sampling points, and blue dots denoting the newly added points.  The errors corresponding to each adaptive step are plotted in Figure \ref{Fig:Test2-2} ($4$ hidden layers are used in the DNN structure).  As we can see, our adaptive sampling strategy accurately captures the singularities and added more sampling points near the interface.  Figure \ref{Fig:Test2-1} compares the exact solution and the DNN approximation.  
In Table \ref{Tab:Test2}, we quantify and compare the performances between the uniform sampling and adaptive sampling strategies.  Again, the adaptive sampling approach uses about $70\%$ fewer points to achieve comparative error. Therefore, we conclude that using adaptive sampling based on the residual of the loss function $\mathcal{L}_{\text{total}}$ is effective for solutions with singularities.

\begin{table}[H]
\caption{Example~\ref{Sect:NumTest2}. $\dfrac{\|u - \mathcal{U_{NN}}\|_\Omega}{\|u\|_\Omega}$. }\label{Tab:Test2}
\centering
\begin{tabular}{c||cccc||cccc}
\hline\hline
&\multicolumn{4}{c||}{4 Layers} &\multicolumn{4}{c}{6 Layers}  \\
Grids & $M_{1}$ & $M_{2}$  & $M_{\Gamma}$ &Error & $M_{1}$ & $M_{2}$  & $M_{\Gamma}$  &Error \\[5pt] \hline
Uniform $10\times 10$ &16 &84 &50 &1.8407e-03 &16 &84 &50 &5.5350e-02\\
Refined Level 1 &201 &96 &52   &6.0430e-04 &164 &133 &53 &4.7112e-03\\
Refined Level 2 &235 &219 &52 &1.3021e-03 &164 &290 &53 &3.7275e-04\\
Refined Level 3 &332 &266 &124 &4.9604e-03 &249 &358 &122 &3.7275e-04\\ \hline
Uniform $50\times 50$  &484 &2016 &160 &7.9539e-03 &484 &2016 &160 &2.3428e-03\\ 
\hline\hline
\end{tabular}
\end{table}

\section{Conclusion}\label{sec:conclusions}
In this paper, we investigate a new deep least-squares method to solve the elliptic interface problem with complicated interface geometries and conditions.  Due to the geometry complexity and/or singularities near the interfaces, classical numerical methods needs special treatment for either meshing technique or modifying the basis functions. Different from previous work, we propose to approximate the solution by deep neural networks and, observing that the solutions might have large jumps in the derivative across the interfaces, we propose to use different DNN structure in each sub-domain.  We then rewrite the interface problem, including the interface and boundary conditions, in the least-squares formulation and the mean squared error loss functions are used on the discrete level so that it can be efficiently trained by the SGD method or its variants.  To capture the singularities, we use the residual error of the loss function as the a posterior error estimator and design an adaptive sampling algorithm. The proposed deep least-squares method is easy to implement and can handle complicated interfaces efficiently. Our numerical experiments show that the proposed deep least-squares method is quite effective for the interface problem and the adaptive sampling strategy improves accuracy while reducing the overall cost for challenging interface problems.




\bibliographystyle{elsarticle-harv}
\bibliography{DLS.bib}

\end{document}